\documentclass[11pt,oneside]{amsart}

\usepackage[backref=page,colorlinks=true]{hyperref}
\usepackage{amsmath,amsthm,amsfonts,amssymb,amscd,url}
\usepackage[capitalize]{cleveref}
\usepackage{graphics,xcolor}
\numberwithin{equation}{section}
\input{xy} \xyoption{all}
\usepackage[utf8]{inputenc}
\usepackage[paper=a4paper, marginpar=2.4cm]{geometry}

\usepackage{enumerate}
\usepackage{hyperref}
\usepackage{epsfig} 
\input{xy} \xyoption{all} 
\usepackage{mathrsfs}

\def\diag {\mathrm{\diag}\, }

\def\Rot{\mathrm{Rot}}
\def~{ }
\def\Ham{\mathrm{Ham}}
\def \Fl {\mathrm{Fl}}
\def \Symp{\mathrm{Symp}}
\def\Re{\mathrm{Re\,}}
\def\Im{\mathrm{Im\,}}

\def\Emb{\mathrm{Emb}}

\def \cA{{\mathscr A}}

\def \cV{{\mathcal V}}

\def \cP{{\mathscr P}}
\def \cT{{\mathcal T}}

\def \cW{{\mathcal W}}

\def\qand{\quad \text{and}\quad}

\def\A{\mathbb A}
\def\B{\mathbb B}
\def\C{\mathbb C}
\def\D{\mathbb D}

\def\M{\mathbb M}
\def\N{\mathbb N}

\def\Q{\mathbb Q}
\def\R{\mathbb R}
\def\V{\mathbb V}
\def\W{\mathbb W}
\def\S{\mathbb S}
\def\T{\mathbb T}
\def\U{\mathbb U}
\def\Z{\mathbb Z}

\def\cU{\mathcal U}

\def\cU{\mathcal U}

\usepackage{mathtools}
\def\trans{\cap\kern-0.7em|\kern0.7em}

\DeclareMathOperator{\leb}{Leb}
\DeclareMathOperator{\Fix}{Fix}

\newtheorem*{theorem*} {Theorem}
\newtheorem*{proposition*} {Proposition}

\newtheorem{proposition}{Proposition}[section]
\newtheorem{theorem}[proposition]{Theorem}

\newtheorem{definition}[proposition] {Definition}
\newtheorem{lemma}[proposition] {Lemma}

\newtheorem{theo}{Theorem}

\newtheorem{fact}[proposition]{Fact}

\newtheorem{corollary}[proposition]{Corollary}
\newtheorem{notation}[proposition]{Notation}

\theoremstyle{definition}
\newtheorem{example}[proposition]{Example}
\newtheorem{remark}[proposition]{Remark}

\newtheorem{problem}[proposition]{Problem}

\begin{document}

\title{Analytic pseudo-rotations II:\\
 a principle for spheres, disks and annuli}
\author[P. Berger]{Pierre Berger}
\address{CNRS IMJ-PRG, Sorbonne Université, Université Paris Cité., partially supported by the ERC project 818737 Emergence of wild differentiable dynamical systems.}

\maketitle

\begin{flushright}
\em To Dmitry Turaev on his 60th birthday. 
\end{flushright}

\begin{abstract} 
We construct analytic surface symplectomorphisms with unstable elliptic fixed points; this solves a problem of Birkhoff (1927). More precisely, we construct analytic symplectomorphisms of the sphere and of the disk which  are transitive, with respectively only  2 and 1 periodic points.
This also solves problems of proposed by Herman (1998), Fayad-Katok (2004) and Fayad-Krikorian (2018).
To establish these results, we introduce a principle that enables  to realize, by an analytic symplectomorphism, properties which are $C^0$-realizable by the approximation by  conjugacy method of Anosov-Katok.
\end{abstract}

\tableofcontents

\section{Main results}
A basic  question in dynamical systems  is the one of stability of a fixed point of an analytic and symplectic, surface dynamics.  
We recall that a \emph{symplectomorphism $F$ of an oriented surface} is an area and orientation preserving diffeomorphism. 
A fixed point $P=F(P)$ is \emph{stable} if all nearby orbits remain close to $P$ for all forward times, i.e., for every neighborhood $V$ of $P$, there exists a neighborhood $U$ of $P$ such that $\bigcup_{n\ge 0} F^n(U)\subset V$.

The problem of the stability of  fixed points of symplectomorphisms is one of the oldest problems in mathematical physics; it goes back to the 18th century. It proposes to characterize the stability of a fixed point  $P$ based on  the properties of the eigenvalues of the differential $D_PF$.  If the modulus of one eigenvalue is not one, then there is an eigenvalue with modulus $>1$  and the fixed point is unstable. 
Examples of unstable fixed points with eigenvalues equal to roots of unity were constructed by Levi-Civita and Cherry \cite{LC01, Ch28}. 
 The fixed point is \emph{elliptic}  if both eigenvalues have  modulus one and are not roots of unity.  

In 1927, Birkhoff \cite[P. 227]{birkhoff1927dynamical} asked wether \emph{all} elliptic periodic orbits of surface analytic symplectomorphisms\footnote{He wrote that an ``outstanding problem'' is wether any elliptic periodic orbit (which is always formally stable) of a Hamiltonian dynamics with two degrees of freedom is stable; using a Poincaré section or conversely an analytic suspension, this problem is equivalent to the stability of  elliptic fixed points of surface analytic symplectomorphisms. He also asked this question for the dynamics given by the  restricted 3-body-problem.} are stable. He presented this problem at many other occasions, for instance  \cite[6th problem]{birkhoff1929einige} as related by the survey \cite{barrow2022george}. 
For \emph{almost all} unit eigenvalues $|\lambda|=1$, the elliptic points are stable by R\"ussmann's theorem \cite{Ru02}. In contrast, for Hamiltonian flow of higher dimension $\ge 6$, 
 Fayad \cite{Fa23}  constructed recently unstable   equilibria of any elliptic eigenvalues. See also \cite{kozlov2023formal}. Yet Birkhoff's problem remained open and of current interest, see \cite[Q. 18 \& Q. 22]{FK18}. 
Our first result solves this problem:
\begin{theo} \label{elliptic instable}There is an analytic  symplectomorphism of the 2-sphere  displaying  an unstable elliptic point. 
\end{theo} 
\begin{figure}[h!] 
\includegraphics[height=7cm]{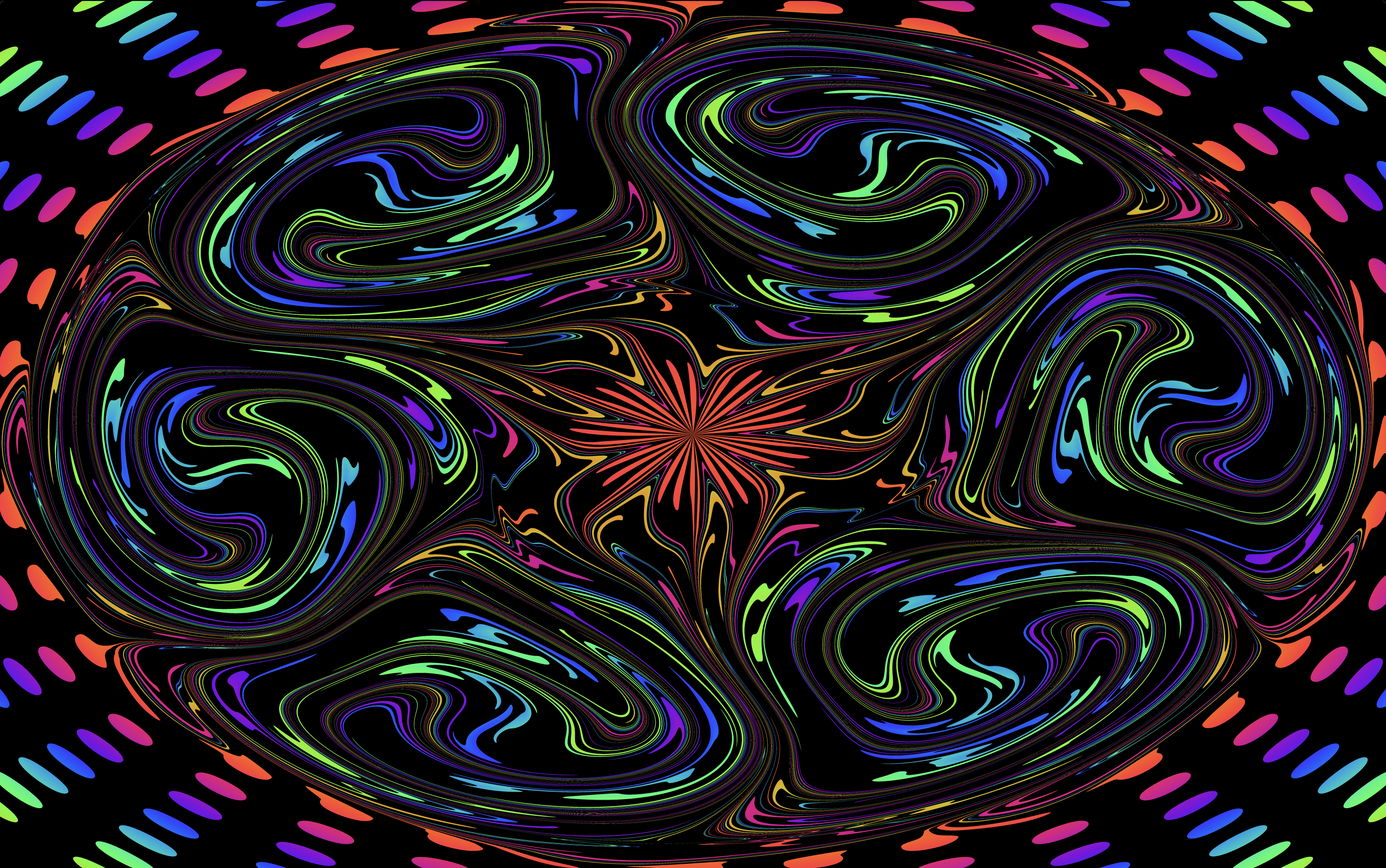}
\caption{First steps of the AbC construction of an unstable elliptic point.} 
\end{figure} 
This result is a consequence of a new principle (see   \cref{h M} \cpageref{h M})  which enables to realize by analytic, surface symplectomorphisms properties which are $C^r$-realizable by the Approximation by Conjugacy (AbC) method of    Anosov-Katok.

The AbC method proved the existence of smooth, transitive, volume preserving dynamics  with a finite number of periodic points on  compact manifolds endowed with an effective action of $\R/\Z$ such as   the cylinder, the disk or the sphere \cite{AK70}. A dynamics is \emph{transitive} if it displays a dense orbit.  A dynamics with a finite number of periodic points is called a \emph{pseudo-rotation}, see  \cite{beguin2006pseudo}. 

Once implemented, the approximation by conjugacy method is extremely efficient in producing dynamics with special properties (e.g. minimal but not ergodic, transitive but not minimal, uniquely ergodic, high local emergence, etc.), see \cite{FK14} for a survey.

However, after the works of Anosov-Katok, it was not clear how this method could be adapted to the analytic setting, see \cite[\textsection 7.1]{FK04}. The  analytic examples were implemented on   
torii \cite[Thm 2.1]{Fu61}, \cite{BK19} or odd spheres \cite{FK14}, but no other manifolds.  In \cite{Be22}, we performed this method in the case of analytic symplectomorphisms of the cylinder. This enabled us to disprove the Birkhoff conjecture \cite{Bi41} by showing the existence of analytic pseudo-rotations of the cylinder or the sphere, which are \emph{not} conjugate to a rotation. The dynamics obtained on the cylinder was transitive.    
Hence, this solved the first half of Herman's problem \cite[Pb 3.1]{He98}, asking whether there exist analytic pseudo-rotations of the cylinder or the sphere with dense orbits. The following complement \cite
{Be22} to answer the second part of Herman's problem:
\begin{theo} \label{Herman problem2}
There exists an analytic symplectomorphism $f$ of the $2$-sphere $\S$  which is transitive and with only two periodic points.  These are elliptic fixed points. 
\end{theo}

Note that \cref{Herman problem2} implies \cref{elliptic instable}. \cref{Herman problem2} is a consequence of the AbC principle (see \cref{h M}). This principle is valid on the cylinder, the sphere and the disk. The latter case enables the construction of analytic, transitive pseudo-rotations of the disk  and hence to solve  problems  presented by Fayad-Katok \cite[Pb. 7.5 and 7.6]{FK04}  and Fayad-Krikorian  \cite[Q. 4]{FK18}:
\begin{theo} \label{Katok problem2}
There exists an analytic symplectomorphism $f$ of the  closed disk which is transitive and has only one periodic point. The fixed point is elliptic.  
\end{theo} 

Theorems A, B and C are the first elementary  applications of the AbC principle, see    \cref{h M} \cpageref{h M}; many other applications are also possible, see for instance \cref{problem C0}. See Problems \ref{pb T2n}, \ref{pb Symp}, \ref{pb vol} and  \ref{finite ergodicity} for some new perspectives offered by the AbC principle.

\thanks{\emph{Thanks:} I am thankful to Olivier Biquard, Penka Georgieva, Julien Grivaux, Pierre Schapira and Dror Varoroline  for many discutions on complex and analytic geometry. I am grateful to  Abed Bounemoura, Bassam Fayad, Pierre-Antoine Guilheuneuf, Vincent Humili\`ere, Raphael Krikorian and  Sobhan Seyfaddini for other  discussions. I am thankful to Curtis McMullen for suggesting me to set up my work in the spirit of the h-principle.
 }

\section{The AbC principle} 
The circle  $\T:= \R/ \Z$ acts  by rotation  $\alpha \in \T\mapsto \Rot_\alpha$ on:
\begin{itemize}
\item  the cylinder $\A:= \T\times [-1,1]$ as $\Rot_\alpha: (\theta,y)\mapsto (\theta+\alpha, y)$,
\item the unit closed disk $\D:= \{ x\in \R^2:  \|x\|\le 1\}$ as the rotation   centered at $0$,
\item  the $2$-sphere $\S:= \{x\in \R^3: \|x\|=1\}$ as the rotation of vertical axis.
\end{itemize}
The surfaces $\A$, $\S$ and  $\D$ are endowed with a canonical symplectic form (of volume 1)  that we all denote by $\Omega$. 
For $1\le r\le \infty$ and $\M\in \{ \A, \D, \S\}$, we denote by  $\Symp^r(\M) $  the space of $C^r$-symplectomorphisms: the subspace of $C^r$-diffeomorphisms  which leave invariant the symplectic form $\Omega$. For $r=0$, we denote by 
$\Symp^0(\M) $  the subspace of homeomorphisms which preserve the volume and the orientation.

We observe that these surfaces are symplectomorphic:
\[ \check \A:= \T\times (-1,1)\;  , \quad \check \D:= \{x\in \D: 0<\|x\|<1\} \qand \check \S:= \{x\in \S: x\neq (0,0,\pm 1)\}\; .\]

Let $\M\in \{ \A, \D, \S\}$.  Let $r\ge 0$ and $\cT\, ^r$ be the  compact-open $C^r$-topology of $\Symp^r(\check \M)$. For instance, we recall that a basis of open neighborhoods of $f\in \Symp^r(\check \A)$ is:
\[ \{ g\in   \Symp^r(\check \A): \| (g-f)| \T\times [-1+\epsilon, 1-\epsilon]\|_{C^s}< \epsilon\} \;  \text{among }\epsilon>0  \text{ and finite integer }s\le r\; .\] 
\begin{definition} 
A \emph{$C^r$-AbC scheme} is a map:
\[ (h,\alpha) \in \Symp^r(\M)\times \Q/\Z \mapsto  (U(h,\alpha), \nu(h,\alpha) )\in   \cT\, ^r \times (0,\infty)\]
such that:
\begin{enumerate}[$a)$]
\item   Each open set $U(h,\alpha)$  contains    the restriction of a map  $\hat h  \in \Symp^r(\M)$ satisfying:  
\[\hat h^{-1} \circ \Rot_\alpha\circ \hat h =  h^{-1} \circ \Rot_\alpha\circ  h\; .\] 
\item   For every sequence $(h_n)_n \in  \Symp^r(\M)^\N$ and $(\alpha_n)_n\in (\Q/\Z )^\N$  satisfying:
\[h_0=id\;,   \quad    \alpha_0=0\; , \]
 and for $n\ge 0$:
\[h_{n+1}^{-1}\circ \Rot_{\alpha_n}\circ h_{n+1} = h_n^{-1}\circ  \Rot_{\alpha_n}\circ h_n \; ,\]\[ \quad 
h_{n+1}|\check \M   \in U(h_n, \alpha_n)\qand 0<|\alpha_{n+1}-\alpha_n|  < \nu(h_{n+1}, \alpha_n),
\] 
then the sequence $f_n:= h_n^{-1} \circ \Rot_{\alpha_n} \circ h_n$ converges to a map $f$ in $\Symp^r(\M)$.
\end{enumerate} 
 The map $f$ is said \emph{constructed from the $C^r$-scheme $(U,\nu)$}.
\end{definition}

\begin{definition}
A property $(\cP)$ on  $\Symp^r(\M)$ is an \emph{AbC realizable $C^r$-property} if:
\begin{enumerate} 
\item  Property $(\cP)$ is invariant by $\Symp^r(\M)$-conjugacy: if $f\in \Symp^r(\M)$ satisfies  $(\cP)$ then so does  $h^{-1}\circ f\circ h$  for any $h\in \Symp^r(\M)$.
\item There exists a  $C^r$-AbC scheme such that every  map  $f\in \Symp^r(\M)$ constructed from it satisfies 
 Property $(\cP)$. 
\end{enumerate}
\end{definition} 
For instance, a classical argument will be used in  \cref{example of ABC} \cpageref{example of ABC} to prove:
\begin{proposition}\label{AbC transitivity ss rot} Being transitive is  an AbC realizable $C^0$-property.
\end{proposition}
Actually, the AbC schemes enable to realize pseudo-rotations. We recall that a \emph{pseudo-rotation} is a symplectomorphism with a finite number of periodic points, which must be   0, 1 or 2  when $\M=\A$, $\D$ or $\S$, and which are moreover all elliptic. 
 We will prove in \cref{section AbC pseudo}:
\begin{proposition}\label{AbC pseudo} 
For every $r\ge 1$, let $(\mathcal  Q)$ be the $C^r$-property of being a pseudo-rotation and let  $(\cP)$ be any  AbC realizable  $C^r$-property. Then the conjunction $ (\cP)\wedge  (\mathcal Q)$ is  an AbC realizable $C^r$-property. \end{proposition}

%
The following is accessible:
\begin{problem}\label{problem C0}
Show that:
\begin{enumerate}
\item Ergodicity is an  AbC realizable $C^0$-property.  
\item  Displaying local emergence of maximal order is an AbC realizable $C^1$-property.
\end{enumerate}
\end{problem} 
As a warm up for the AbC principle stated below, we can prove:
\begin{proposition} \label{CRC=Cinfty}
For every $0\le s\le r\le \infty$, an  AbC realizable  $C^s$-property  is an AbC realizable $C^{r}$-property. In particular an AbC realizable $C^0$-property is an AbC realizable $C^\infty$-property.
\end{proposition} 
\begin{proof} Let $(\cP)$ be a property realized by a $C^s$-AbC scheme  $(U,\nu)$. 
 Also $\Symp^r(\M)$ is a Fréchet space for a distance $d$. 
For $\ell= s$ or $r$, let $\Symp^\ell_c(\check \M)$ be the dense  subspace of $\Symp^\ell(\check \M)$ formed by symplectomorphisms which coincide with $id $ outside of a compact subset of $\check \M$.   We recall that  $\Symp^r_c(\check \M)$ is dense  in $\Symp^s_c(\check \M) $ by    \cite{Oh06,Si07}  and \cite{Ze77}. Using a finite covering, for every $q\ge 1$, the subspace of  $\Symp^r_c(\check \M)$ formed by maps commuting  with $\Rot_{1/q}$ is dense  in the subspace of $\Symp^s_c(\check \M) $  formed by maps commuting  with $\Rot_{1/q}$.   Hence, for every 
$(h,\alpha)\cap \Symp^r(\M)\times \Q/\Z$, the open set $U(h,\alpha)\cap \Symp^r ( \check \M)$ contains the restriction of a map $\hat h\in  \Symp^r _c(  \M)$ such that:
\[\hat h\circ  \Rot_{\alpha}\circ \hat h=h ^{-1}\circ \Rot_{\alpha}\circ h  \; .\]
Now let us consider a   sequences  $(h_n)_n \in  \Symp^r(\M)^\N$ and $(\alpha_n)_n\in (\Q/\Z )^\N$  satisfying:
\[ h_0=id\;  , \; \alpha_0=0\; , \]
and  for $n\ge 0$:
\[h_{n+1}^{-1}\circ \Rot_{\alpha_n}\circ h_{n+1}= h_n^{-1}\circ  \Rot_{\alpha_n}\circ h_n  \; ,\]
\[ \quad 
h_{n+1}|\check \M   \in U(h_n, \alpha_n)\cap \Symp^r(\check \M)\qand 0<|\alpha_{n+1}-\alpha_n|  < \nu(h_{n+1}, \alpha_n)\; .
\] 
Up to taking $\alpha_{n+1}$ closer to $\alpha_n$, we can assume that 
$f_{n+1}=h_{n+1}^{-1}\circ \Rot_{\alpha_{n+1}}\circ h_{n+1} $ is $C^r$-$2^{-n}$-close to $f_{n}=h_{n}^{-1}\circ \Rot_{\alpha_{n}}\circ h_{n} $.
This corresponds to assume $|\alpha_{n+1}-\alpha_n| $ smaller than a certain function  $\tilde \nu(h_{n+1}, \alpha_n)$.  Then observe that $(f_n)_n$ converges to a map $f\in \Symp^r (\M)$. We accomplished to show that $(U\cap  \Symp^r(\check \M) , \tilde \nu)$ is  a $C^r$-scheme  which realizes property $(\cP)$. 
 \end{proof}
 
 From the three latter propositions, we deduce:
\begin{proposition}\label{AbC transitivity} Being a transitive pseudo-rotation is  an AbC realizable $C^1$-property.
\end{proposition}

Here is the main theorem of this work:
\begin{theo}[AbC principle] \label{h M}
For every $0\le r\le \infty$,  any AbC realizable $C^{r}$-property $(\cP)$ is satisfied by a certain analytic  symplectomorphism $f\in \Symp^\omega(\M)$.
\end{theo} 
Note that this theorem together with \cref{AbC transitivity} in the case $\M= \S$ implies  immediately \cref{Herman problem2}, while in the case of $\M=\D$, it implies immediately \cref{Katok problem2}. 

It is natural to ask on which other manifolds the AbC principle holds true. 

By \cite{BT22},  every symplectomorphism of $\T^{2n}$ can be $C^\infty$-approximated by an entire symplectomorphism.  From this we easily deduce that a certain AbC principle holds true in $\Symp^\omega(\T^{2n})$.
\begin{problem}\label{pb T2n} State properly  and show the AbC principle on  $\Symp^\omega(\T^{2n})$.
\end{problem} 

 Some AbC constructions exist also on any toric\footnote{manifold $(M^{2n},\Omega)$ displaying an effective action of $\T^n$.}   symplectic manifolds, see \cite{LRS22}.  This together with the techniques presented in this work should be useful to find other examples of  symplectic manifolds satisfying the  AbC principle:
\begin{problem}\label{pb Symp}
For which symplectic manifolds an AbC principle holds true? 
\end{problem}  
Another way to generalize the AbC principle is within the category of volume preserving analytic maps, see \cite{FK14} for the case of odd spheres. 
\begin{problem} \label{pb vol}
For which other manifolds endowed with a volume form an  AbC principle holds true?
%
\end{problem} 
 
It seems that the AbC principle as stated above does not encompass  finite ergodicity. Yet this should be reachable:
\begin{problem}\label{finite ergodicity}
State and prove an AbC principle which enables to prove that   finite ergodicity is realizable by an analytic symplectomorphism of the sphere, the disk or the cylinder.

\emph{Hint:} Instead of approximating some  $h|\check \M$ in the AbC scheme,   given any pair of open sets  $M_1\subset \check \M$ and $M_2\Subset \M\setminus M_1$, one 
should approximate  $h|M_1$ by a map which is close to the identity on $M_2$.  
A combination of the techniques presented here enables to prove such an approximation  result. 
\end{problem}

\section{Structure of the proof of the main theorem} 

The main difficulty for implementing the AbC scheme in the analytic setting is the lake of analytic symplectic mappings  which 
commute with a rational rotation and extending holomorphically to a large complexification of the surface.
 For instance, an easy consequence of Picard's theorem is that any analytic symplectic map of $\A=\R/\Z\times [-1,1]$ which extends to $\C/\Z\times \C$ must be a shear of the form $(\theta,y)\mapsto (\theta+g(y), y)$.  

To overcome  this rigidity, the idea is to deform the analytic structure of these surfaces. 

Basically, we will extend the conjugacies to complex differentiable maps which are close to be holomorphic and close to leave invariant the holomorphic extension of the symplectic form. These conjugacies  pull balk the canonical complex structure on a complexification of the surface to form a sequence  that we will show converging to a  (possibly different)  complex structure. For this complex structure, the realization of the scheme will extend holomophically as well as the symplectic form. This will imply  that the realization is analytic and symplectic for a (possibly different) structure. This latter sketch of proof will be more detailed  in \cref{sec sketch}.  In the next  \cref{sec:Rigidity of the real analytic symplectic structure}, we recall that up to isomorphism, there is a unique  analytic and symplectic structure.

\subsection{Rigidity of the real analytic symplectic structure}\label{sec:Rigidity of the real analytic symplectic structure}
We recall that a \emph{smooth structure  of symplectic surface}  $(\M,\Omega)$ is a maximal $C^\infty$-atlas $\cA^\infty=(\phi_\alpha)_\alpha$  such that each chart  $\phi_\alpha : U\subset \M\to V\subset \R^2$ pulls back  the symplectic form $dx\wedge dy$ to $\Omega|U$.
A \emph{real analytic structure of symplectic surface}  $(\M,\Omega)$ is a maximal  atlas $\cA^\omega=(\phi_\beta)_\beta$  such that 
the coordinate changes are real analytic and 
each chart  $\phi_\alpha : U\subset \M\to V\subset \R^2$ pulls back the symplectic form $dx\wedge dy$ to $\Omega|U$. 
We say that the structure  $\cA^\omega$ is \emph{compatible} with $\cA^\infty$ if $\cA^\omega\subset \cA^\infty$.

 \begin{proof}[Proof of \cref{h M} with $\M=\S$] By \cref{CRC=Cinfty}, it suffices to show the case $r=\infty$. Then the main difficulty will be to prove:
\begin{theorem}\label{Main1}    For any AbC realizable $C^{\infty}$-property  $(\cP)$ on $\Symp^\infty(\S)$, there exists a real analytic symplectic structure $\cA^{\omega'}$ on $(\S,\Omega)$  which is compatible with the canonical  $C^\infty$-symplectic structure  and  for which  $(\cP)$ is satisfied by an analytic  symplectomorphism $f$.\end{theorem}
Yet the inclusion $\S\subset \R^3\subset \C^3$ endows  $(\S,\Omega)$ with a canonical real analytic symplectic structure $\cA^\omega$ that we would like to recover. To this end  we  infer:
\begin{theorem}[\cite{KL00}] \label{rigidity1}
Any symplectic manifold possesses a unique (up to isomorphism) real analytic symplectic structure.\end{theorem}
Hence there exists an analytic symplectomorphism $\Psi: (\S, \cA^\omega) \to (\S,\cA^{\omega'})$. Then observe that $\Psi^{-1} \circ f  \circ \Psi$ is a $C^\omega$-symplectomorphism which satisfies  $(\cP)$. 
\end{proof} 

To prove Theorems \ref{h M} in the case  $\M= \A$ or $\D$, we would like to use a counterpart of  \cref{rigidity1} for manifolds with boundary. To this end, there is the following theorem on analytic symplectic structures which are equivariant by an analytic  action of a compact Lie group $K$ (with its natural real analytic structure):
\begin{theorem}[{\cite{KL00} \textsection4}] \label{rigidity2}  Any symplectic K-manifold has a unique real analytic $K$-structure.\end{theorem}
We can use this theorem because $\M$ has a canonical structure of symplectic $\Z_2$-manifold. Indeed,  we can glue two copies of $\M$ along its boundary to form a boundaryless surface $\tilde \M$ so that $\M$ is the quotient of $\tilde \M$ by a reflection. 

For instance, $(\A, \Omega)$ is canonically identified to the quotient of the torus  $\hat \A:= (\T\times \R/4\Z, \Omega)$ by  the involution: 
\[\tau: (\theta,y)\in  \T\times \R/4\Z \mapsto (\theta, -y+2).\]
Likewise, $(\D, \Omega)$ is canonically identified to the quotient of a symplectic sphere $(\tilde \D, \Omega) $ by an involution denoted also by $\tau$, as we will see in \cref{tilde S}. 

To apply   \cref{rigidity2} we shall  consider the real analytic $\Z_2$-structure on $\M\in \{\D, \A\}$. A \emph{real analytic $\Z_2$-structure of symplectic surface}  $(\M,\Omega)$ is a maximal  real analytic symplectic atlas $\cA^\omega=(\phi_\beta)_\beta$  such that each chart  $\phi_\alpha: U_\alpha\to \R\times \R_+ $ can be lifted to an open set  $\tilde U_\alpha\subset \tilde  \M$ to form a $\Z_2$-equivariant chart $\tilde \phi_\alpha$ of   $\tilde \M$.  This means that with
$\tau_0(x_1,x_2):=  (x_1, -x_2)$, the following diagram commutes:
\[  
\begin{array}{lcccr}
&\tilde U_\alpha & \stackrel{\tau}\longrightarrow & \tilde U_\alpha&\\
\tilde \phi_\alpha&\downarrow & & \downarrow &\tilde \phi_\alpha\\
&\R^2 & \stackrel{\tau_o}\longrightarrow & \R^2& 
 \end{array}  \; . \]

\begin{proof}[Proof of \cref{h M} with $\M=\A$ or $\D$] By \cref{CRC=Cinfty}, it suffices to show the case $r=\infty$.
The main difficulty will be to prove:
\begin{theorem}\label{Main2}   
 For any AbC realizable $C^{\infty}$-property $(\cP)$ on $\Symp^\infty(\M)$, there exists a $\Z_2$-real analytic symplectic  structure $\cA^{\omega'}$ on $ \M$  which is compatible with the  canonical  $C^\infty$-structure  and  for which   $(\cP)$ is satisfied by    a $\Z_2$-analytic  symplectomorphism $f$.
\end{theorem}
The inclusions $\A\subset \R/\Z\times \R/4\Z \subset \C/\Z\times \C/4\Z$ and $\D\subset \R^2 \subset \C^2$ endow $  \M$ with a canonical real analytic symplectic $\Z_2$-structure $\cA^\omega$ that we can recover. Indeed, by \cref{rigidity2},  there exists an analytic $\Z_2$-symplectomorphism $\Psi: (   \M , \cA^\omega) \to ( \M ,\cA^{\omega'})$. Then 
 $\Psi^{-1} \circ f  \circ \Psi$ is a $C^\omega$-symplectomorphism of $\M$  which satisfies  property~$(\cP)$. 
 \end{proof}

\subsection{Sketch of  proof of the main theorem modulo  deformation}\label{sec sketch}
The proof of  \cref{Main1} and  \cref{Main2} will occupy all the remaining of this manuscript.  

In \cref{Space}, we define and introduce the notations of the involved  spaces:  the disk, the  sphere, the cylinder, and their complexifications. For instance,   the cylinder is $\A:= \R/\Z\times [-1,1]$ and   the sphere is $\S=\{ x\in \R^3: \sum_i x_i^2=1\}$.
We will work with the following complexifications of the sphere: \[\S_\infty:= \{ z\in \C^3: \sum _i z_i^2=1\}\qand \S_\rho:= \{ z\in \S_\infty: \sum |z_i|^2<\rho\}\quad \text{for } \rho>1\; .
\]

The idea is to use the main approximation theorem of \cite{Be22}    which enables to approximate any symplectic map of $\check \A:= \R/\Z \times (-1,1)$ into $\R/\Z\times \R$ by  an entire automorphism of  $\R/\Z \times \R$, whose restriction to a large complex neighborhood of  $\R/\Z \times \{-1,1\}$ in $\C/\Z\times \C$    is close to the identify.  A natural way to import this result to the sphere is to use the axial projection:
\[  (x_1,x_2,x_3)\in \S \setminus \{(0,0,,\pm1) \} \mapsto (\tfrac1{2\pi} \arg (x_1+i x_2) , x_3)\in \R/\Z\times (-1,1)\; .\] 
The axial projection projection is well known to be symplectic.  In \cref{Space}, we will show that  it extends to a biholomorphism $\pi_\S$ between the two following spaces:
\[\pi_\S:  \check \S_\infty:=\{ (z_1,z_2,z_3) \in \S_\infty: |\Re(z_3 )|<1\}\stackrel{\sim}{\rightarrow} \check \A_\infty:=\{  (\theta , y)\in \C/\Z\times \C  : |\Re(y )|<1\}\; .\]

In \cref{structure}, we will give a sufficient condition for a smooth symplectomorphism $f$ to be analytic for a compatible structure.   This will be given by 
\cref{Def struture analytico symplectic} in term of  (integrable almost) complex structure $J$. For instance, in the case of the sphere,  we  recall that a complex structure   $J$ on $\S_\rho$
 is an automorphism of its real tangent bundle $T^\R\S_\rho$ whose square is $J^2=-id$ and such that   there exists a $C^\infty $-atlas of $\S_\rho$ formed by \emph{$J$-holomorphic charts}
 $\phi_\alpha$. This means that the charts satisfy $D\phi_\alpha \circ J= i D\phi_\alpha$.  
Then \cref{Def struture analytico symplectic} will imply the following:
\begin{proposition*} Let $\rho>\rho'>1$. Let $J$ be a  complex structure on $\S_{\rho}$ which anti-commutes with $(z,u)\in T\S_{\rho}\mapsto (\bar z,\bar u)$
and such that the symplectic form $\Omega$ of $\S$ extends $J$-holomorphically to $\S_\rho$.  Then $J$ defines  a real analytic symplectic structure on $\S$ which is compatible with the $C^\infty$-structure.    Furthermore if a symplectic map $f$ of $\S$ has a 
 $J$-holomorphic extension   $ \S_{\rho'} \to \S_{\rho}$, then $f$ is   real analytic and symplectic for this structure. \end{proposition*}

Let us explain how useful will be this proposition.  First let us  endow the space of complex structures on   $\S_\rho$   with  the topology induced by the $C^\infty$-topology on the space  smooth  sections   of $(T^\R \S_\rho)^*\otimes T^\R \S_\rho$. We will see in \cref{cpx sr closed}   that this space is complete (i.e. the integrability condition   is $C^\infty$-closed). 
Also we will see in \cref{CS holo} that the existence of a holomorphic extension of the symplectic form is also a closed property in the $C^\infty$-topologies on complex 2-forms and  complex structures. Hence to prove  \cref{Main1} , it will suffice to realize the AbC scheme using conjugacies which extend to semi-conjugacies $(H_n)_n$ from $\check \S_\rho:=\check \S_\infty\cap \S_\rho$ onto their image in $\C/\Z\times \C$, such that 
 each $H_n$ satisfies $H_n(\bar z)= \overline{H_n(z)}$ and each $H_n$ pulls back the canonical complex structures $J_o$ of $\C/\Z\times \C$ and its holomorphic 2-form $\Omega_o:=d\theta\wedge dy$ in way that $(H^*_nJ_o)_n 
 $ and $(H^*\Omega_o)_n$ converge in the $C^\infty$-topology. Then the limit of $(H^*_nJ_o)_n$ is  a complex structure $J$   and  the limit of 
 $(H^*\Omega_o)_n$ is 
 a $J$-holomorphic form $\Omega$. By the above proposition, for the structure induced by $J$, the restriction of    $f= \lim H_n^{-1} \circ \Rot_{\alpha_n} \circ H_n$ to $\S$ is an analytic symplectomorphism. From this we will deduce  \cref{Main1}.  
%
%
%
%
%
%

We will start \cref{sec:approxa} by stating a  development  of the main approximation result of   \cite{Be22}, as \cref{approxa}. This Theorem involves the following complexifications of the cylinders $\A:= \T\times [-1,1]$ and $\check \A:= \T\times (-1,1)$, indexed by $\infty\ge R>1$:
\begin{multline*}
 \A_\infty:= \C/\Z\times \C \; , \quad \A_R := \{ (\theta, y)\in \A_\infty : |\exp( i \theta)|^2+ |\exp( -i \theta)|^2+ |y|^2\le 3R\}\\
 \qand \check \A_R := \{ (\theta,y)\in \A_R: |\Re(y|<1\}\; .
\end{multline*}
A  rough version of  \cref{approxa} is:

  \begin{theorem*}  Let $R>1$ and let $\tilde h$ be a symplectomorphism of $\check \A $ which commutes with a rational rotation $\Rot_\alpha$. 
  Then there exists a smooth diffeomorphism   $h:  \A_R \hookrightarrow \A_\infty$  onto its image  such that:
\begin{enumerate}[($M_1$)]\setcounter{enumi}{-1}
\item $h$ commutes with $\Rot_\alpha$,
\item $h( \A)=  \A$,    $ \det Dh| \A =1$ and  $h|\check \A$ is $C^\infty$-close to $\tilde h$,
\item the map  $ h$ coincides with the identity on the complement of $\check \A_R$,
\item the pullbacks  $h^* J_o $  and $h^*\Omega_o$ are  close to $J_o$ and $\Omega_o$, 
\item  the map $h$ commutes with $\sigma:(\theta,y)\mapsto (\bar \theta, \bar y)$.
\end{enumerate} 
\end{theorem*}
From this  theorem,  we will deduce \cref{Main2} in the cylinder case in \cref{AbC realization on deformed real analytic symplectic cylinders}. Indeed, we will show that  given an AbC scheme and compact neighborhoods $\A_{\rho'}\Subset \A_{\rho}$ of $\A$ in $\C/\Z\times \C$, the later result enables to construct a sequence of maps $h_n:\A_{\rho} \to \A_\infty$ and of rational numbers $\alpha_n$ such that:
\begin{itemize}
\item  The diffeomorphism $h_n$ leaves invariant $\A$, the restriction $h_n|\A$ is symplectic and $(h_n|\A)_n$ follows the AbC scheme with $(\alpha_n)_n$. 
\item The sequence $(h_n^*J_o)_n$ converges to  a complex structure $J$ on $\A_{\rho}$.
\item The sequence $(h_n^*\Omega_o)_n$ converges to  a $J$-holomorphic  form  on $\A_{\rho}$.
\item The diffeomorphism  $f_n:= h_n^{-1} \circ \Rot_{\alpha_n} \circ h_n| \A_{\rho'}$ is well defined and $(f_n)_n$  converges to a smooth map~$f$.
\item The diffeomorphism  $h_n$ commutes with $\sigma$.
\end{itemize}
From this and the results of \cref{structure},  we will deduce that  $J$ defines  a real analytic symplectic structure on a neighborhood of $\A$ in $\R/\Z\times \R$, which is compatible with the $C^\infty$-structure and which is preserved by  $f|\A_{\rho'} \cap \R/\Z\times \R$.  
Then the main  \cref{Main2} in the cylinder case is obtained by 
developing this construction to obtain a $\Z_2$-structure.

\begin{figure}[h!] 
\includegraphics[height=7cm]{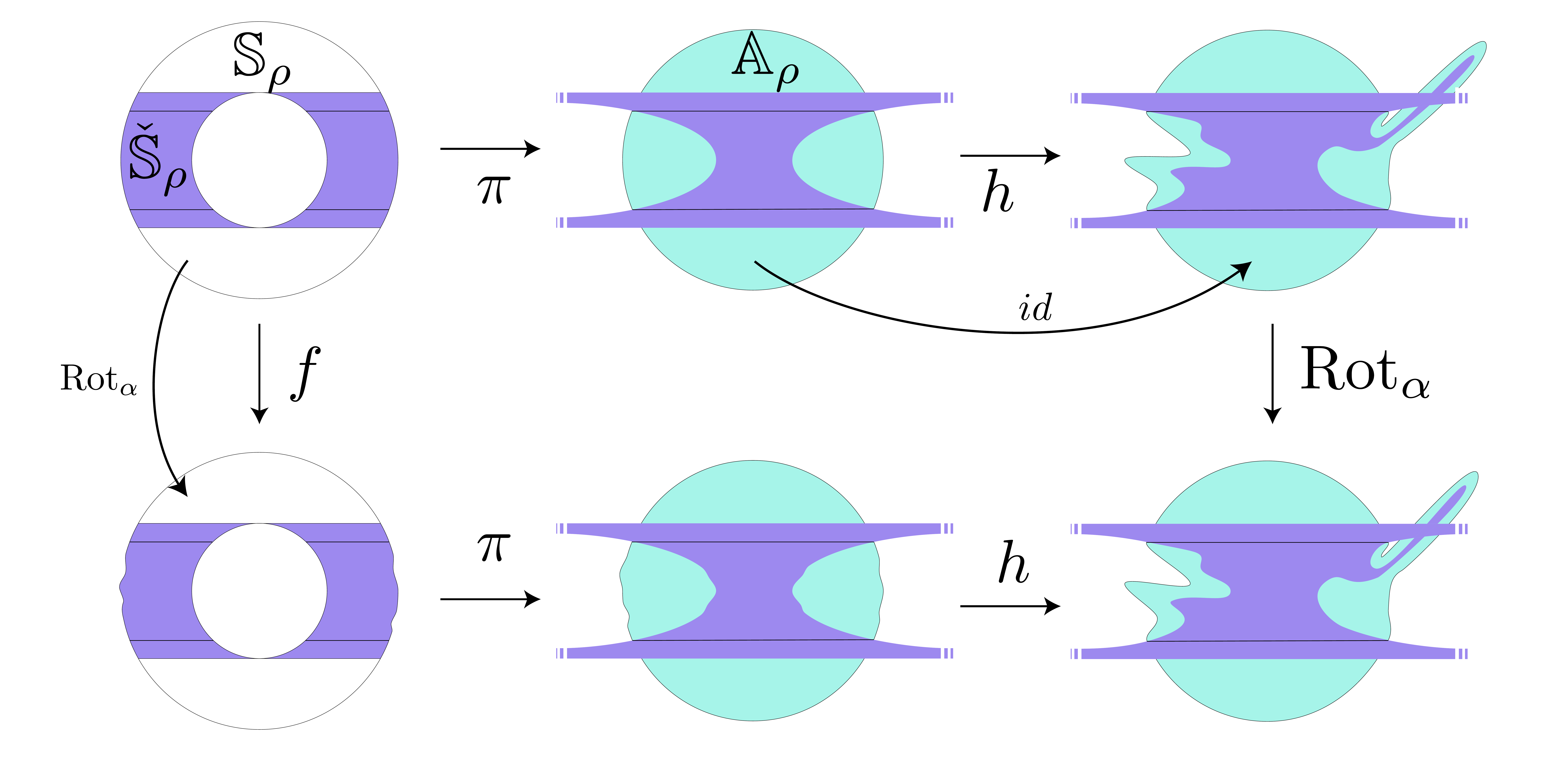}
\caption{AbC realization for deformed real analytic symplectic spheres.}\label{sketch_Sphere}
\end{figure}
In \cref{AbC realization for deformed real analytic symplectic spheres and disks}, we will show how a development of this construction proves the sphere and disk cases of \cref{Main1}. Let us sketch the sphere case.  
Consider a sequence $(h_n)_n$ as above.  Then the map $ h_n\circ \pi$ is well defined on a certain $\check \S_r$ intersected with $\pi^{-1}(\A_\rho)$.  Also as $h_n$ coincides with the identity on a neighborhood of the complement of $\check \A_\rho$, for suitable values of $r$ and $\rho$, see \cref{sketch_Sphere},  we can extend $h_n$ by the identity on $\pi(\check \S_r)\setminus \A_\rho$.  Then $H_n:= h_n\circ \pi$ is well defined on $\check \S_\rho$ and coincides with $\pi_\S$ nearby  $\S_r \setminus \check \S_r$.  
So the pullback $H_n^*J_o$ can be extended by the   canonical  complex structure on  $\S_r\setminus \check \S_r$.  This defines a complex structure  $H_n^\#J_o$ on $\S_r$. One  defines similarly the form $H_n^\#\Omega_o$. Assuming the convergences of  $(h_n^*\Omega)_n$ and $(h_n^*J)_n$ are ``sufficiently fast",  we will obtain that $H_n^\#\Omega$ and $H_n^\#J$ also converge. However, as depicted by  \cref{sketch_Sphere},   the maps $H_n:  \check \S_r\to \C/\Z\times \C$ are in general non-injective  local diffeomorphisms. Yet we will manage to extend Anosov-Katok's tricks to this setting. This will define a sequence of maps $f_n: \S_{r'}\to \S_r$ whose restrictions to $\check \S_{r'}$ are  semi-conjugate to  rational rotations  $\Rot_{\alpha_n}$ via $H_n$ and which coincide to $\Rot_{\alpha_n}$ abroad. These maps are  $H_n^\#J_o$-holomorphic. We will show that for a ``sufficiently fast'' convergence of  $(\alpha_n)_n$, the sequence $(f_n)_n$ converges to a diffeomorphism $f$.  Then $f|\S$ realizes the scheme AbC, and by the same structure theorems of \cref{structure}, the map $f|\S$ leaves invariant an  analytic  symplectic structure on $\S$, as sought by   \cref{Main1}. The disk case is obtained similarly. 

We finish  in \cref{proof approx} by proving   \cref{approxa}. In order to do so we will recall the main approximation theorem of \cite{Be22}. A rather easy consequence  of which   is \cref{coro approx analytic} which states roughly   the following.   Given a compact neighborhood $\A_R$ of $\check \A := \R/\Z\times (-1,1)$ with smooth boundary, any symplectic map $\tilde h $ from $\check \A$ into $\R/\Z\times \R$ can be approximated by the restriction of an entire automorphism $h_\bullet$ of  $\R/\Z \times \R$, so that $ h_\bullet$ is $C^\infty$-close to the identity on a neighborhood of  $ \{(\theta,y)\in \A_R: |\Re(y)|\ge 1\}$.   Also if $\tilde h$ commutes with a rational rotation $\Rot_{\alpha}(\theta,y)=(\theta+\alpha, y)$, then $h_\bullet$ can be chosen commuting with $\Rot_\alpha$.   Then  the idea is to  use a bump function to glue $h_\bullet$ to the identity  of $ \{(\theta,y)\in \A_R: |\Re(y)|\ge 1\}$ to define a map $h$ which is $C^\infty$-close to $h_\bullet$ and so close to be holomorphic and leaving invariant $\Omega_o$. 
Finally we perturb it so that its real trace leaves $\A$ and the symplectic form invariant.

\section{Spaces and their notations}\label{Space}
In this section, we specify the expressions of the cylinder, the sphere and the disk. We also define their canonical complexifications and their $\Z_2$-structure if any. 
Moreover we set up the symplectic  axial projection $\pi_\S$ and the polar symplectic coordinates $\pi_\D$ together with their holomorphic extensions.

\subsection{Complexifications of the cylinder $\A$}
The cylinder and its interior are:
\[\A:= \T\times [-1,1]\qand \check \A:= \T\times (-1,1)\; .\] 
For $\infty\ge\rho>1$, we have the following complexifications:
\[ \A_\infty:= \C/\Z\times \C \; , \quad 
 \A_\rho := \{ (\theta, y)\in \A_\infty: |\exp( i \theta)|^2+ |\exp( -i \theta)|^2+ |y|^2\le 3\rho\}
\]
\[\qand  \check \A_\rho :=\{ (\theta,y)\in \A_\rho: |\Re y| <1\}\; .\]
Note that $\A$ is the quotient by the involution $\tau : (\theta, y)\mapsto (\theta,2-y)$ of 
\[\tilde \A:=  \T\times \R/4\Z\; .\]
The involution $\tau$ extends holomorphically to $\A_\infty$.   Let $\tilde \A_\rho:= p(\A_\rho)\cup \tau \circ p(\A_\rho)$ where $p$ is the projection $\A_\infty\to \C/\Z\times \C/4\Z$. 
Observe that $\tilde \A_\rho$ is left invariant by $\tau$. 
 
 We endow $\A$ with the symplectic form $\Omega:= \tfrac12 d\theta\wedge dy$, which extends holomorphically to $\A_\infty$ as a holomorphic form denoted by $\Omega_o$.  Note that $\tau $ is symplectic.  

\subsection{Complexifications of the  sphere  $\S$}
We recall that:
\[\S:= \{x\in \R^3: \sum x_i^2=1\}\qand \check \S:= \{x\in \S: |x_3|<1\}\; .\]
We denote again by $\Omega$ the canonical symplectic form on $\S$ with volume $1$.  The twice punctured sphere $\check \S$ is symplectomorphic to the cylinder $\check \A$  via the axial projection:
\[ x\in \check \S \mapsto (\tfrac1{2\pi} \arg (x_1+ix_2) , x_3)\in \check \A\; .\]
 For $\rho>1$, we have the following complexifications:
\[ \S_\infty:= \{ z\in \C^3: \sum z_i^2=1\}\; , \quad 
 \S_\rho := \{ z\in \S_\infty: \sum |z_i|^2\le \rho\}\; ,\]
\[\qand \check \S_\rho := \{ z\in \S_\rho: |\Re z_3| <1\}\; . 
\]

We extend analytically the square root   to $\C\setminus \R_-$. Then for any $z\in \check \S_\infty$,
as $|\Re z_3|<1$, the real part of $z_3^2$ is smaller than 1 and so we have $z_1^2+z_2^2=1-z_3^2\in \C\setminus \R_-$. Thus, we can consider $\sqrt{z_1^2+z_2^2}$. 
Likewise we recall that $\arg$ is a function from the unit circle of $\C$ onto $\R/2\pi\Z$ which extends analytically to $\C^*$ as the inverse  $-i  \log : \C^*\to \C/2\pi\Z$  of the  map $z\mapsto \exp( i  z)$. This enables to extend the axial projection by   the following   bi-holomorphism:
\[ \pi_\S: \check \S_\infty \to
\check \A_\infty =\{ (\theta,y)\in \A_\infty: |\Re y| <1\} \]
\[ z\mapsto \left ( \frac1{2 i \pi}\log  \left( \frac { z_1+i z_2}{\sqrt {z_1^2 +z_2^2}}\right) , z_3\right) \]
with inverse:
\[(\theta,y)\in \check \A_\infty \mapsto \left(\sqrt{1-y^2} \cdot (\cos (2\pi \theta) ,\sin (2\pi \theta)  ) ,y\right)\; .\]
We recall that $\pi_\S|\check \S$ is a symplectomorphism onto $\check \A$.  As the symplectic form of $\check \A$ extends holomorphically to $\check \A_\infty$, the symplectic form  of $\S$ extends holomorphically to $\check \S_\infty$. By symmetry, it extends also holomorphically to any image of  $\check \S_\infty$ by an element of $SO_3$. As $SO_3(\check \S_\infty)= \S_\infty$ and since $\S_\infty$ is simply connected\footnote{$\S_\infty$ is diffeomorphic to $T\S$.}, the symplectic form $\Omega$ of $\S$ extends holomorphically to $ \S_\infty$ as a holomorphic $2$-form denoted also by $\Omega_o$. 

\subsection{Complexifications of the disk $\D$} \label{tilde S} We recall that:
\[ \D :=  \{ x\in \R^2 : \sum x_i ^2\le 1\}\qand \check  \D :=  \{ x\in \R^2 : 0<\sum x_i ^2< 1\}\; .
\]
We denote again $\Omega:=\tfrac1\pi dx_1\wedge dx_2$ the canonical symplectic form on $\D$. 
The  punctured disk $\check \D$ is symplectomorphic to the cylinder $\check \A$  via the symplectic polar coordinates:
\[ x\in \R^2\setminus \{0\}  \mapsto \left(\frac1{2 \pi}\arg (x_1+i x_2) , 2(x_1^2+x_2^2)-1\right)\; .\]

For $\rho>1$, we have the following complexifications:
\[ 
 \B_\rho := \{ z\in \C^2: \sum |z_i|^2\le \rho\} \qand  \check \D_\rho := \{ z\in \B_\rho: 0<\Re (z_1^2+ z_2^2) <1\}\; , 
\]
on which the symplectic form extends holomorphically as $\Omega_o:=\tfrac1\pi dz_1\wedge dz_2$. 
Likewise,   we extend the symplectic polar coordinate by the following  bi-holomorphism:
\[ \pi_\D: \{ z\in \C^2: \Re(z_1^2+z_2^2)>0\}\to \{ (\theta,y)\in \A_\infty: \Re y>-1\} \]
\[ z\mapsto \left (\frac1{2 i \pi}\log  \left( \frac { z_1+i z_2}{\sqrt {z_1^2 +z_2^2}}\right) , 2 (z_1^2+z_2^2)-1\right) \]
We notice that $\pi_\D$ sends   $\check \D$ onto $\check \A$, and sends  
$\check \D_\infty $ to $\check \A_\infty$. Let us describe a $\Z_2$-equivariant structure on $\D$. Let $\U$ be the neighborhood of $\partial \D$ defined by:
\[\U:= 
\{x\in \R^2: 0<x_1^2+x_2^2<2\}\; .\]
Note that $\\U$ is sent by $\pi_\D$ onto $\T\times (-1,3)$ which  is left invariant by the symplectic involution  $ (\theta, y)\in \T\times \R \mapsto (\theta, 2-y)$.  The latter map  is conjugate via $\pi_\D|\U$ to the following  involution of $\U$:
\[ \psi:= x\in \U \mapsto  \sqrt{2-\|x\|^2}\cdot \frac{ x}{\|x\|}=\sqrt{\frac 2{x_1^2+x_2^2} -1}\cdot x \; .
\]
As $\pi_\D|\U$ is symplectic, the map $ \psi$  is a symplectic involution of $\U$; it leaves invariant $\partial \D$.

Then  the real analytic symplectic structure on $\D\cup \U$ induces canonically a   real analytic symplectic surface structure  (see  \cite[Cor. 2.4]{BeCoex23}) on the following quotient:
\[ \tilde \D:=  (\D\cup \U)\times \Z_2  /\sim,\]
where $(x,i)\sim(x', i')$ if $(x,i)=(x', i')$ or  if $x\in \U$ and  $(x',i') = (\psi(x),i+1)$. Note that the manifold $\tilde \D$ is a symplectic sphere. 
Also  $\D$ is the quotient of $\tilde \D$ by the involution $\tau: [x, i ] \in  \tilde \D \mapsto [x, i+1]\in \tilde \D $.

Let us now propose a canonical complexification of $\tilde \D$.  To this end, using again  the analytic extension of $\sqrt{\cdot } $ on $\C\setminus \R_-$, 
we  observe that  $\psi $ extends analytically as:
\[  \Psi: z\in \{ (z_1,z_2)\in \C^2: 0<\Re(z_1^2+z_2^2)<2\}  \mapsto z\cdot \frac {\sqrt{2- \sum z_i ^2}}{\sqrt{\sum z_i ^2}} 
 \; .\] 
 The map $\Psi$ is still an involution since it is conjugate by $\pi_\D$ to the involution:
 \[(\theta, y)\in \C/\Z \times \{ 0< \Re(y)<2\}\mapsto (\theta,2-y)\; .\]

Let 
\[\U_{\rho}:= \{ z\in \B_\rho:  -1<\Re(z_1^2+z_2^2) <3\}\; .\]
We define:
\[ \D_\rho := \{ z\in \B_\rho: \Re(z_1^2+z_2^2)<2\} \; .\]
We notice that the following is a complex extension of $\tilde \D$:
\[ \tilde \D_\rho :=  \D_\rho\times \Z_2 /\sim \quad \text{with }
 (z, i)\sim (\Psi(z), i+1)\text{ when }z\in  \U_{\rho} \; .\] 
We can extend the   involution $\tau$ as $[z, j ] \in  \tilde \D_\rho \mapsto [z, j+1]\in \tilde \D_\rho $. It leaves invariant the holomorphic form $\Omega_o$ on  $\tilde  \D_\rho$, which  equals to $ \tfrac1\pi dz_1\wedge dz_2$ in  canonical coordinates.

\section{Structures and deformations}\label{structure}
In this section we define complexifications of  a smooth symplectic manifolds and set a complete topology on the space of such complexifications. We finally 
give useful criteria for surface maps  to be real analytic and symplectic,  for the  structure induced by such complexifications. 

\subsection{$C^\infty$-Topology} 
Let $M$ and $N$ be manifolds (possibly with boundary). The topology  on $ C^\infty(M,N)$  is given by the union of the compact-open $C^r$-topologies among $r\ge 1$. 
In other words, a sequence  $f_n\in C^\infty(M,N)$  converges to $f\in C^\infty(M,N)$ iff for every compact domain $W\subset M$ with smooth boundary and $r\ge 1$,  the sequence of restrictions $f_n|W$ converges uniformly to $f|W$ in the $C^r$-topology.  

\subsection{Complex structures and their topologies} 
Let $\V$ be a complex manifold.   We denote by $T\V$ the complex tangent bundle of $\V$ and by $T^\R\V$ its structure of real vector bundle.   
The multiplication by $i$ on the complex vector bundle structure defines an automorphism $J_o$ of $T^\R\V$ such that $J_o^{2}=-id$.  

\begin{notation} 
The manifolds $\C^k$, $\A_\infty$, $\tilde \A_\infty$, $\D_\infty$, $\tilde \D_\infty$ and $\S_\infty$ have canonical structures of complex manifolds that we all denote by $J_o$.  
\end{notation}

More generally, a $(C^\infty)$ smooth section $J$ of $T^\R\V^*\otimes T^\R\V$  such that $J^{2}=-id$ is called an  \emph{almost complex structure}. 
An almost complex structure $J$ on a manifold is called \emph{integrable} if
there is a $C^\infty$-atlas $(\phi_\alpha)_\alpha$ of $\V$ formed by charts $\phi_\alpha: U_\alpha \to \R^{2n}\approx \C^n$ such that $D\phi_\alpha \circ J= i D\phi_\alpha$. Then observe that the coordinate changes of the atlas are holomomorphic. 

Hence an integrable  complex structure defines a complex structure on $\V$ and vice-versa. Therefore  we  identify  the   complex structures  to   the  integrable almost complex structures. 

Let $\W\subset \V$ be a compact domain with smooth boundary. We denote by $T^\R\W$ the restriction $T^\R\V|\W$. We still denote by $J_o$  the complex structure induced on $\W$.   
 \begin{example}\label{ex cpx strcture} 
For every  smooth embedding  $h$  from $\W$ into $\V$, the following 
\[ h^*J_o:= Dh^{-1}\circ J_o\circ Dh \]
 is a complex structure on  $\W$. 
 \end{example} 

We endow    the space of  complex structures on $\W$ with the   topology  induced by    the space of smooth sections of  $T^\R\W^*\otimes T^\R\W$.  The following is classical:

\begin{theorem} \label{cpx sr closed}
Let $(J_n)_n$ be a sequence of complex structures on $\W$ which converges to an almost complex structure $J$ in the $C^\infty$-topology. Then $J$ is integrable: it is  a complex structure on $\W$.   
\end{theorem}\begin{proof}
Given an   almost complex structure $J$, the  Nijenhuis tensor $N_{J}$ is the tensor which  to a pair $(X,Y)$ of smooth sections of  $T^\R\V$  associates the following section:
\[
N_{J'}(X,Y) = [X,Y]+J([JX,Y]+[X,JY])-[JX,JY]\; .
\] 
The Newlander–Nirenberg theorem states that an almost complex structure $J$ is  integrable iff the tensor $N_J=0$.
The nullity of the Nijenhuis tensor can be checked on finitely many vector fields. Hence the integrability condition is closed in the $C^\infty$-topology.
\end{proof} 
 
 \subsection{Complexification of real analytic structure} 
Let $\M$ be a real analytic manifold.  A \emph{complexification} of $\M$ is a complex manifold $( \M_\bullet,J)$ such that $\M_\bullet\supset \M$ and such that for every $z\in \M$, there exists a $J$-holomorphic chart $\psi_\alpha$ of a neighborhood   $V_\alpha$ of $z$ in $\M_\bullet$  whose restriction   to $U_\alpha=\M\cap V_\alpha$ is a real analytic chart of $\M$. It is easy to show that every real  analytic manifold displays a complexification. In this subsection, we   give  a sufficient condition for a complex manifold  $( \M_\bullet,J)$ containing a smooth real manifold $\M$  to be a complexification of a certain   real analytic structure on $\M$.

A \emph{real structure} on a complex manifold $(\M_\bullet,J)$  is an involution $\sigma: \M_\bullet\to \M_\bullet$ which is \emph{antiholomorphic}:
\[ \sigma^2=id \qand  \sigma^*  J= -J \; .\]
We will endow the spaces $\M_\bullet= \A_\rho, \D_\rho, \S_\rho,\C^k, \tilde \A_\rho, \tilde \D_\rho ...$  with the canonical real structure:
\[ z=(z_j)_j\in \V\mapsto \bar z= (\bar z_j)_j\in \V\; .\]

\begin{proposition}
\label{pre Def struture analytico symplectic}
Let $\M\in \{\A, \D, \S,\tilde \A, \tilde \D\}$, $\rho>1$ and  let $J$ be a complex structure on $\M_\rho$  such that $\sigma^*  J=-J$.  Then there exists a real analytic structure   on $\M$ which is compatible with the $C^\infty$-structure of $\M$ and for which  $(\M_\rho, J)$ is a  complexification.  

Moreover if $\M=  \tilde \A$ or $ \tilde \D$ and  $D\tau \circ J= J\circ D\tau $, then  the latter structure   contains a $\tau$-equivariant atlas which projects to a real analytic  $\Z_2$-structure on $\A$ or $\D$.
\end{proposition}
\begin{proof} 
For every $z\in \M\subset \M_\rho$, there is a distinguished neighborhood $V_\alpha$ of $z$ and a chart $\psi_\alpha: V_\alpha \to \C^2$ which is $J$-holomorphic. Up to composing $\psi_\alpha$ with a translation, we can assume that $\psi_\alpha(z)= 0$. 
Now take a basis $(e_1, e_2)$ of $ T_z \M$. Up to composing $\psi_\alpha$ with a linear map, we can assume that $e_1$ and $e_2$ are sent to $(1,0)$ and $(0,1)$.  Then observe that the following map sends also the basis  $(e_1,e_2)$  to $((1,0),(0,1))$:
\[\phi_\alpha:=  \frac12(\psi_\alpha + \sigma\circ \psi_\alpha \circ \sigma)|U_\alpha\quad  \text{with } U_\alpha\text{ a neighborhood of $z$ in }V_\alpha\cap \M\; .\]
Hence for $U_\alpha$ sufficiently small,  $\phi_\alpha$ is a local diffeomorphism onto its image in $\R^2$. Observe that $  \phi_\alpha$ commutes with $\sigma$. This implies that $\phi_\alpha$ sends $U_\alpha$ into $\R^2$.  We  observe that such maps $(\phi_\alpha )_\alpha$  form a real analytic  atlas of $\M$ and which is compatible with the $C^\infty$-structure of $\M$.

If $\M=  \tilde \A$ or $\tilde \D$ and $z\in \Fix\, \tau $, then we chose a basis $(u,v)$ of $T_z \M$ such that $D_z\tau(v)=-v$ and $u\in T_z \Fix\,  \tau$. In particular, $D_z\tau(u)=  u$. 
Up to a composition by a real linear map, we can assume that $D_z \phi_\alpha (u)= (1,0)$  and $D_z \phi_\alpha (v)= (0,1)$. Then the following 
\[\tilde  \phi_\alpha:=  \frac12( \phi_\alpha + \tau_o\circ  \phi_\alpha \circ \tau) \]
 sends also $z$ to $0$,  sends  also $u$ and $v$ to respectively $(1,0)$ and $(0,1)$ and  still  commutes with $\sigma$ (because $\tau$ does). Moreover 
$\tilde  \phi_\alpha\circ \tau= \tau_o\circ \tilde  \phi_\alpha$.  So the latter kinds of maps  descend to form a $\Z_2$-real analytic atlas on $\M$. 
\end{proof} 
By the latter proposition, the following is consistent:
\begin{definition} A complex structure $J $ on $\M_\rho$ such that $ \sigma^* J  = - J$ is called a  \emph{complexification of a real analytic structure} on $\M$.\end{definition}

\begin{example}\label{ex cpx strcture2} Let $h$ be a $C^\infty$-embedding of $\M_\rho$ into $\M_\infty$ which commutes with $\sigma$. Then $(\M_\rho, h^*J_o)$ is a complexification of $\M$. 
Indeed we saw that  $h^*J_o$ is a complex structure in \cref{ex cpx strcture}. Also it satisfies  $\sigma^*h^*J_o=h^*\sigma^*J_o=h^*(-J_o)=-h ^* J_o$, so we can apply the latter proposition.
  \end{example}
  \begin{remark} If $h$ is a $C^\infty$-embedding of $\M_\rho$ into $\M_\infty$ which commutes with $\sigma$, then $h(\M)\subset \M$ because $\M=\Fix\sigma$.
  \end{remark}
\subsection{$J$-holomorphic extensions of  symplectic forms} 
In this subsection we give criteria on the symplectic form $\Omega$ to be analytic for the structure induced by a complexification.

 \begin{definition} Let $\M$ be a real surface endowed with a smooth symplectic form $\Omega$ and let  $(\M_\bullet, J)$ be a complexification of  $\M$. 
 We say that   $\Omega$ extends \emph{$J$-holomorphically to $\M_\bullet$}, 
 if there is a \emph{$J$-holomorphic   extension} $\Omega_\bullet$ of $\Omega$ on $\M_\bullet$, i.e. $\Omega_\bullet|T\M=\Omega$  and  seen in   $J$-holomorphic charts, the form $\Omega_\bullet$ equals to:
 \[ \psi(z_1,z_2)\cdot dz_1\wedge dz_2,\]
  for a holomorphic function $\psi$ on an open subset of $\C^2$.
 \end{definition}
 Let  $\M\in \{\A, \D, \S,\tilde \A, \tilde \D\}$. We recall that $\Omega_o$ denotes the (canonical) $J_o$-holomorphic extension on $\M_\infty$ of the symplectic form $\Omega$ on $\M$.  Let $\rho>0$.
\begin{example}\label{example extension symp} Let $h$ be a $C^\infty$-embedding of $\M_\rho$ into $\M_\infty$ which commutes with $\sigma$ and whose restriction to $\M$ is symplectic (it leaves invariant $\Omega$).  Then 
$\Omega$ extends $h^*J_o$-holomorphically to $\M_\rho$.  Indeed $h^*\Omega_o$ is an extension of  $\Omega$. Also an atlas of  $(M_\rho,h^*J_o)$ if formed by charts $\psi_\alpha\circ h$, where $\psi_\alpha$ is a chart of $(M_\rho, J_o)$. Then note  
that seen in this chart, $h^*\Omega_o$ is $(\psi_\alpha\circ h)_* h^*\Omega_o= (\psi_\alpha)_*\Omega_o$ which is indeed holomorphic. 
  \end{example} 
A \emph{complex differential 2-form} on a manifold  is a  differential $2$-form  which is permitted to have complex coefficients.  We endow the space of complex differential  2-forms  on a compact manifold with the $C^\infty$-topology. The following is classical:
\begin{theorem} \label{CS holo} 
 Let $(\M,\Omega)$ be a real symplectic surface. Let  $(\M_\rho, J_n)_n$ be a sequence of complexifications of  $\M$  which converges to a complexification $(\M_\rho, J_\infty)$, i.e. $J_n\to J_\infty$. 
 Assume that   $\Omega$ has a  $J_n$-holomorphic   extension $\Omega_n$  on $\M_\rho$ for every $n$ and that $(\Omega_n)_n$   converges.  Then $\Omega$ extends $J_\infty$-holomorphically to $\M_\rho$.
\end{theorem}\begin{proof} Note that $\Omega_\infty=\lim_\infty \Omega_n$ is a smooth extension of $\Omega$. It suffices to show that $\Omega_\infty$  is $J_\infty$-holomorphic. 
To this end we prove below:\begin{lemma}\label{CS holo2} 
A   complex differential  2-form  $\Omega_\bullet$ on $\M_\rho$ is holomorphic iff it is closed and $\Omega_\bullet(J u,v)=i\Omega_\bullet (u,v)$ for any pair of vectors $(u,v)$.
 \end{lemma}
Then  for every $n$, we have $\Omega_n(J_n\cdot, \cdot)=i\Omega_n (\cdot, \cdot)$, and so at the limit we have $\Omega_\infty(J_\infty\cdot, \cdot)=i\Omega_\infty(\cdot, \cdot)$. Using again the lemma, we obtain that  $\Omega_\infty$ is $J_\infty$-holomorphic.  
\end{proof}
 \begin{proof}[Proof of \cref{CS holo2}]
 Clearly, if $\Omega_\bullet$ is $J$-holomorphic then $\Omega_\bullet$ is closed and $\Omega_\bullet(J \cdot , \cdot)=i\Omega_\bullet $.   
 Now if $\Omega_\bullet(J \cdot , \cdot)=i\Omega_\bullet $, then a simple computation shows that  in a $J$-holomorphic  chart the form cannot have component in $ dz_1\wedge d \bar z_2$, $ \bar dz_1\wedge d  z_2$ or $d \bar z_1\wedge d\bar  z_2$. Hence there exists a smooth function $g$ such that $\Omega_\bullet (z)  =  g  (z)\cdot  dz_1\wedge dz_2$. If moreover $\Omega_\bullet$ is closed, then:
%
 \[ 0=d\Omega_\bullet   =  \frac\partial{\partial \bar z_1} g   \cdot d\bar z_1\wedge  dz_1\wedge dz_2+ \frac\partial{\partial \bar z_2} g   \cdot d\bar z_2\wedge  dz_1\wedge dz_2\; ,\]
and so $\bar \partial g=0$. This implies that $g$ is holomorphic. 
 \end{proof}
\subsection{$J$-holomorphic extensions of symplectomorphisms}
Let $(\M_\rho ,J)$ be a complexification of a  real manifold $\M$.   Let $\M_{\rho'}\subset \M_\rho$ be a smaller complexification of $\M$.  
In this subsection we give a useful criterion for a map $f: \M_{\rho'}\to \M_\rho$ to define a real analytic symplectomorphism of $\M$  for the structure induced by $J$.

\begin{definition} 
A map $f\in C^\infty(\M_{\rho'}, \M_\rho)$ is a \emph{$J$-holomorphic extension of a real map} if $f(\M)\subset \M$   and  $f$  is  \emph{$J$-holomorphic}:  $Df\circ J= J\circ Df$.
\end{definition}

\begin{proposition}  \label{Def struture analytico symplectic}Let $\rho>\rho'$. 
Let $(\M_\rho,J)$ be a complexification of $\M\in  \{ \S, \tilde \A, \tilde \D\}$. 
Assume that $\Omega$ extends $J$-holomorphically to $\M_\rho$.  
Let $f\in C^\infty(\M_{\rho'}, \M_\rho)$ be a  $J$-holomorphic extension of a real symplectomorphism of $\M$.  Then $f|\M$ is a   real analytic symplectomorphism for a $C^\infty$-compatible symplectic structure. 

 Moreover, if $\M=  \tilde \A$ or $\tilde \D$,  assume that $\tau^*  J= J $ and that  $f$ commutes with $\tau$.  Then $f|\M$  is the lifting of an analytic $\Z_2$-symplectomorphism on $\A$ or $\D$.
\end{proposition}
\begin{proof}  In the real analytic chart $(U_\alpha, \phi_\alpha)$ given by the real analytic structure induced by $(M_\rho,J)$, there is  a real analytic function $\psi_\alpha$ such that  $\Omega$ has the form
\begin{equation}\label{Psialpha}  x=(x_1,x_2)\in U_\alpha \subset \R^2\mapsto \psi_\alpha(x) \cdot dx_1\wedge dx_2,\end{equation} 
 Then, we can compose each chart $\phi_\alpha$ by the inverse of  $g_\alpha:=(x_1,x_2)\mapsto (x_1, \int_{0}^{x_2} \psi_\alpha(x_1, t)dt)$ in order to obtain a chart which sends $\Omega$ to 
 $  dx_1\wedge dx_2$. Clearly $f|\M$ is a real analytic symplectomorphism for this structure. 
 
 In the case of  $\tau^*J=J$, we can assume that $\tau_o\circ \phi_\alpha= \phi_\alpha\circ \tau$   by \cref{pre Def struture analytico symplectic}. Then the function $\psi_\alpha$ of \cref{Psialpha} satisfies $\psi_\alpha \circ \tau_o = \psi_\alpha$, and so $g_\alpha$ commutes with $\tau_o$. Hence $\tilde \phi_\alpha:= g_\alpha^{-1} \circ \phi_\alpha$ is    $\Z_2$-equivariant. Clearly $f|\M$ is a real analytic $\Z_2$-symplectomorphism for this structure. 
  \end{proof}
  \begin{corollary} \label{Def struture analytico symplectic2}Let  $\M\in \{ \A,  \D\}$ and $\rho>\rho'$. 
Let $(\M_\rho,J)$ be a complexification of $\M$. Assume that $\Omega$ extends $J$-holomorphically to $M_\rho$. 

Let $f\in C^\infty(\M_{\rho'}, \M_\rho)$ be a  $J$-holomorphic extension of a real symplectomorphism of $\M$. Assume that $J$  
coincides with $J_o$ on $\M_\rho \setminus \check \M_\rho$.  Assume that $f$  coincides with a rotation $\Rot_\alpha$ on  $\M_{\rho'} \setminus \check \M_{\rho'}$. 

Then  $f|\M$ is  a  real analytic $\Z_2$-symplectomorphism for a  $C^\infty$-compatible structure. 
\end{corollary}
 \begin{proof} Case $\M= \A$.  We recall that  $p$ is the projection $\A_\infty\to \C/\Z\times \C/4\Z$ and that   $\tilde \A_\rho:= p(\A_\rho)\cup \tau \circ p(\A_\rho)$ is endowed with the $\Z_2$-structure given by the involution  $\tau : (\theta, y)\mapsto (\theta,2-y)$.
We notice that:
\[\tilde \A_\rho:= p(\check \A_\rho)\sqcup \tau \circ p(\check \A_\rho)\sqcup p(\partial \check \A_\rho)\; .\] 
We define $\tilde J$ as equal to $p _* J$ on  $p(\check \A_\rho)$, to $(\tau \circ p)_* J$ on  $(\tau \circ p)(\check \A_\rho)$ and to $J_o$   on  $p(\partial \check \A_\rho)$. As 
$p _* J$ extend smoothly by $J_o=\tau^* J_o$ on the complement of $p(\check \A_\rho)$, it holds that  $\tilde J$ is smooth and integrable\footnote{Its Nijenhuis tensor is zero since it is zero for $J$ and $J_o$}. Observe also that $\tilde J$ is $\Z_2$-equivariant and satisfies $\sigma^*\tilde J= -\tilde J$. 
As $J= J_o$ on $  \A_\rho \setminus \check \A_\rho$, the $  J$-holomorphic extension $\Omega_\bullet $ of $\Omega$ coincides with $\Omega_o$ on $\A_\rho \setminus \check \A_\rho$.  Hence we can define $\tilde \Omega$ as equal to $p _* \Omega_\bullet$ on  $p(\check \A_\rho)$, to $(\tau \circ p)_*  \Omega_\bullet$ on  $(\tau \circ p)(\check \A_\rho)$ and to $\Omega_o$   on  $p(\partial \check \A_\rho)$.  Observe  that $\tilde  \Omega$ is smooth and $\tau$-equivariant. 
Let $q$ be the inverse of $p|\check \A_{\rho'}$.  Let $\tilde f$ be the map of  $C^\infty(\tilde \A_{\rho'}, \tilde \A_\rho)$ equal to $p \circ f\circ q$ on $p(\check \A_{\rho'})$,  to $\tau \circ p \circ f\circ q\circ \tau $ on $\tau \circ p(\check \A_{\rho'})$ and  to
$\Rot_\alpha$ on  $p(\partial \check \A_{\rho'})$.  Observe that $\tilde f$ is smooth, $\tilde J$-analytic,  $\tau $-equivariant and has a $\Z_2$-quotient which coincides with $f$ on $\A$.  

Case $\M= \D$. In \cref{tilde S}, we defined:
\[\D_\rho =\{z\in \B_\rho: \Re(z_1^2+z_2^2)<2\}\; ,\quad \check \D_\rho =\{z\in \B_\rho: 0<\Re(z_1^2+z_2^2)<1\}\; ,\]
\[\quad \U_\rho =\{z\in \B_\rho: 0<\Re(z_1^2+z_2^2)<2\}\qand \Psi (\check \D_\rho)
=  \{z\in \C^2: 1<\Re(z_1^2+z_2^2)<2\}\; .\] 
Also we defined: 
   \[ \tilde \D_\rho :=  \D_\rho\times \Z_2 /\sim \quad \text{with } (z, i)\sim (\Psi(z), i+1)\text{ when } z\in \U_\rho\; ,\] 
   and the involution  $\tau:  [z, j ] \in  \tilde \D_\rho \mapsto [z, j+1]\in \tilde \D_\rho $.  We now define:
   \[ \hat J: (z,j)\in \D_\rho\times \Z_2\mapsto \left\{ \begin{array} {cl} 
J(z) & \text{if }  \Re(z_1^2+z_2^2)<1\\
J_0(z) & \text{if }  \Re(z_1^2+z_2^2)=1\\
\Psi^*J(z) & \text{if }  \Re(z_1^2+z_2^2)>1
\end{array}
\right. \; .\]
  As $\Psi ^*J_o=J_o$ coincides with $J$ on the complement of $\check \D_\rho$, the structure $\hat J$ is well defined, smooth and integrable. As the structure $\hat J| \U_\rho\times \Z_2$ is equivariant by $ (z, j)\mapsto  (\Psi(z), j+1)$, it projects to a structure $\tilde J$ on $\tilde \D_\rho$. We notice that $\tilde J$  is  $\Z_2$-equivariant and projects to $J$ on $cl(\check \D_\rho)$.  We proceed similarly to construct and so to show the existence of  $\Z_2$-equivariant 
 $ \tilde J$-holomorphic extensions of $\Omega$ and $f|\D$ to respectively 
  $\tilde \D_\rho$ and  $\tilde \D_{\rho'}$.   \end{proof} 
\section{AbC realizations on deformed real analytic surfaces} \label{sec:approxa} 
 In this section we are going to achieve the proof of the main theorems by proving    \cref{Main1,Main2}, stating the AbC principle for  $C^\infty$-compatible  real analytic symplectic structures. In order to prove this,  in \textsection \ref{sec6.1}, we will introduce the new approximation \cref{approxa} on the cylinder; it will be proved at \textsection  \ref{sec6.4}. In \textsection  \ref{sec6.2} , this approximation theorem will enable to implement for any $\rho>1$,  the AbC scheme on the cylinder with conjugacies 
 displaying    smooth extensions to  $\A_\rho$ which are close to be $J_o$-holomorphic and close to leave invariant the canonical 2-form $\Omega_o$. Using the geometrical propositions of the previous section, we will obtain that the sequence of pullbacks of the structures $J_o$ and $\Omega_o$ by these conjugacies will converge to a complex structure $J$ such that $\Omega_o$ extends $J$-holomorphically to $\A_\rho$. Also the realization of the AbC scheme will have a $J$-holomorphic extension to $\A_\rho$. This will suffice to obtain the cylinder case using \cref{AbC_cyl_1}. In \textsection  \ref{sec6.3},
we will prove the   sphere case and the disk case: $\M\in \{ \S,\D\}$. We will face several difficulties. To overcome them,  we will  pull back the structure  $J_o$  by maps of the form $g_n\circ \pi_\M$ where $g_n$ is a local diffeomorphism of the cylinder constructed inductively using the aforementioned approximation theorem. The proof will need  more cares   because the AbC scheme will be implemented by transiting between two different spaces $\M_\rho$ and $\A_R$ and moreover  because the AbC scheme will be implemented with possibly not injective local diffeomorphisms.   

\subsection{Approximation modulo structural   deformation}
\label{sec6.1} In this subsection, we state the new approximation \cref{approxa} on the cylinder.
We recall that $\A_\infty=\C/\Z\times \C$ is endowed with a canonical complex structure $J_o$ and a canonical holomorphic form $\Omega_o:=\tfrac12  d\theta\wedge dy$. In \cref{Space}, we also  defined   the compact neighborhood $\A_R$ of $\A$ in $\A_\infty$ for every $R>1$. 

\begin{definition} \label{def approx modulo}
We say that $h_o\in \Symp^\infty(\A)$ is \emph{approximable modulo structural deformation}  if for every 
 $C^\infty$-neighborhood $\cU\subset \Symp^\infty(\check \A)$  of $ h_o|\check \A$, then there exists $\eta>0$ such that  the following property holds true. 

For every $R>1$, 
for any $C^\infty$-neighborhoods 
  $\cW_J$ of $J_o| \A_{R} $ and   $\cW_\Omega$ of $\Omega_o| \A_{R} $,  then  there exists a smooth diffeomorphism   $h:  \A_{R} \hookrightarrow \A_\infty$  onto its image  such that:
\begin{enumerate}[($M_1$)] 
\item $h(\A)= \A$,    $ \det Dh|\A =1$ and  $h|\check \A\in \cU$,  
\item the map  $ h$ is supported by   $ \{ (\theta, y)\in \A_{R}: |\Re (y)|<1-\eta\}$,
\item the pullbacks  $h^* J_o $  and $h^*\Omega_o$ are  in  $\cW_J$ and $\cW_\Omega$,
\item  the map  $ h$  commutes with $\sigma$.
\end{enumerate} 
\end{definition} 
This notion will be useful because,  by Examples \ref{ex cpx strcture}, \ref{example extension symp}  and \cref{Def struture analytico symplectic}  then  any map of the form $h^{-1} \circ \Rot_\alpha\circ h|\A$ is analytic for the 
real  analytic symplectic structure on $(\A,\Omega)$  induced by  $(h^*J_o,h^*\Omega_o)$.  
By \cref{Def struture analytico symplectic2,}, it is also true for the induced $\Z_2$-structure. 
 
  Here is the new approximation theorem:
  \begin{theorem} 
 \label{approxa} 
Every $h_o\in \Symp^\infty(\A)$ is approximable modulo structural deformation.
\end{theorem}This theorem  will be shown  in \cref{proof of approxa} by using the main approximation theorem of \cite{Be22}.  Using a lifting $(\theta,y)\mapsto (q\cdot \theta, y)$, an immediate consequence  is:
 \begin{corollary}If $h_o$ commutes with a certain  $\Rot_{p/q}$,  then the approximations $h$ modulo structural deformation can be chosen commuting with  $\Rot_{p/q}$.
 \end{corollary}

\subsection{AbC realizations on deformed real analytic symplectic cylinders}\label{AbC realization on deformed real analytic symplectic cylinders}\label{sec6.2}
 
In this subsection, we prove  \cref{Main2}  in the case $\M=\A$  which states that every $C^\infty$-AbC scheme 
$(U,\nu)$ on $(\A,\Omega)$  is realizable by an analytic symplectomorphism for a  $C^\infty$-compatible real analytic symplectic structure.
This theorem is a consequence of the following proposition shown below.

\begin{proposition}\label{AbC_cyl_1} Let $ \rho>\rho'>1$.  
Then  for every $C^\infty$-AbC scheme $(U,\nu)$,  there is 
 a complexification $(\A_\rho, J)$ on which  $\Omega$ extends $J$-holomorphically and  there is   a $J$-holomorphic  embedding  $f :   \A_{\rho'} \hookrightarrow \A_{\rho} $ which  extends a real symplectomorphism of $(\A,\Omega)$ realizing the   AbC scheme $(U,\nu)$.
 
 Moreover on the complement of $ \check \A_\infty$, the map $f$ equals  a rotation and   $J=J_o$. 
\end{proposition} 
\begin{proof}[Proof of \cref{Main2}  in the case $\M=\A$] It   is 
  a direct consequence of  \cref{AbC_cyl_1} and \cref{Def struture analytico symplectic2}. \end{proof}
 \begin{proof}[Proof of  \cref{AbC_cyl_1}]
 The space of smooth embeddings $\Emb^\infty (\A_{\rho}, \A_\infty)$ is    metrizable for a complete distance $d$.
The space  of complex differentiable  2-forms  on  $\A_{\rho}$  is  metrizable for a complete distance denoted also by $d$.
Likewise the space of smooth linear automophisms of $T^\R\A_{\rho}$ is    metrizable for a complete distance again by~$d$.

We are going to construct $f$ as a limit of maps $f_n$ given by an induction.   More precisely, we will show by induction on $n\ge 0$, the existence  of $\alpha_n\in \Q/\Z$, $h_n\in \Emb^\infty(\A_{\rho}, \A_\infty)$ with support in $\check \A_{\rho}$  and  $f_n\in \Emb^\infty(\A_{\rho'}, \A_\infty)$ such that:
\[ \alpha_0=0\; ,\quad  h_0= id|\A_\rho\qand  f_0= id|\A_{\rho'}\; .\]
And for $n\ge 1$,  $\alpha_n$, $h_n$ and $f_n$ satisfy:
\begin{enumerate}[$(I_n)$]
\item  $h_n(\A)=\A$ and the restriction  $ h_n|\A$ is a symplectomorphism which  respects the AbC scheme $(U,\nu)$  with $\alpha_n$:
\[h_n^{-1}\circ  \Rot_{\alpha_{n-1}}\circ h_n =h_{n-1}^{-1}\circ  \Rot_{\alpha_{n-1}}\circ h_{n-1}\quad \text{on  }\A  \; , \]
\[\;   h_n|\check \A\in U( h_{n-1}| \A,\alpha_{n-1}) \qand  |\alpha_n -\alpha_{n-1}|<\nu(  h_{n}| \A,\alpha_{n-1})\; .
\]
\item The map $ \Rot_{\alpha_n} \circ h_n$ sends $\A_{\rho'}$ into the interior of $h_n(\A_{\rho})$ and so 
the composition $f_n:= h_n^{-1} \circ \Rot_{\alpha_n} \circ h_n| \A_{\rho'}$  is well defined. Moreover  $d(f_n,f_{n-1})<2^{-n}$.
\medskip
\item  The complex structure $J_n:=h_n^*  J_o  $ and the form $\Omega_n:= h_{n}^*\Omega_o$ satisfy:
\[d(J_n,J_{n-1})<2^{-n}\qand d(\Omega_n,\Omega_{n-1})<2^{-n}\; .\]
\medskip
 \item  The maps $h_n$ and $\sigma$ commute:  $h_n\circ \sigma= \sigma\circ h_n$. 
\medskip\item The support of $h_n$ is included in $\check \A_\rho$. 

\end{enumerate} 
 We will prove the following below:
 \begin{lemma}\label{induction} If  the induction hypothesis $(I_n\cdots V_n)$ is satisfied at step $n\ge0$, then it is satisfied at step $n+1$. 
\end{lemma} 
 
Now let us show that the induction hypothesis implies the proposition. 
By  $(II_n)$,  the sequence $(f_n)_n$ converges to a certain $f\in \Emb^\infty(\A_{\rho'}, \A_{\rho})$.  By $(I_n)$, the map $f|\A$ is a symplectomorphism of $\A$ which  realizes the AbC scheme.  

 By $(III_n)$,  the sequence $(J_n)_n$ converges to  an almost complex structure $J$. By \cref{cpx sr closed}, $J$ is a complex structure.  Also by $(I_n)$ and \cref{example extension symp}, the complex form  $\Omega_n$ is a $ J_n$-holomorphic extension of $\Omega$. 
By $(III_n)$ and  \cref{CS holo}, the sequence   $(\Omega_n)_n$ converges to  a $J$-holomorphic extension of $\Omega$. Also 
\[f_n^*J_n= (h_n^{-1} \circ \Rot_{\alpha_n}\circ h_n )^* h_n^*J_o= h_n^* \Rot_{\alpha_n}^*J_o=J_n.\]
 Hence at the limit $f^*J=J$.  Consequently, the map $f$ is a $J$-holomorphic extension of  the symplectomorphism  $f|\A$. 
  Finally by $(V_n)$,  on the complement of $ \check \A_\infty$, the map $f$ equals  a rotation and   $J=J_o$. 
\end{proof} 
 
\begin{proof}[Proof of \cref{induction}]
Let $R>1$ be sufficiently large so that the interior of $ \A_R$ contains  $ h_n(\A_{\rho})$:
\begin{equation} \label{inclusion W 1}  \A_{R}\Supset  h_n(\A_{\rho})\; .\end{equation}
 Note that $\A_{R}$ is $\Rot_{\alpha_n}$-invariant $ \A_{R}= \Rot^{-1}_{\alpha_n} (\A_{R})$ and so:
 \begin{equation} \label{inclusion W 2}
   \A_{R}\Supset  \Rot_{\alpha_n} \circ h_n(\A_{\rho})\; .
 \end{equation} 

By $(III_n)$, there is a $C^\infty$-neighborhood $  \cW_J$ of $J_o|  \A_{R}$  such that for every $J\in  \cW_J$, the  structure $h_n^*  J$ is at distance $\le 2^{-n-1}$ from  $J_n $. 
Also there is a $C^\infty $-neighborhood $ \cW_\Omega$ of $\Omega_o|  \A_{R} $ such that for every $\Omega'\in \cW_\Omega$, the  form  $h_n^* \Omega'$   is  at distance $\le 2^{-n-1}$ from  $\Omega_n $. 

As the scheme $(U,\nu)$ is $C^\infty$-realizable, there exists $\tilde h\in \Symp^\infty(\A)$  which commutes with $\Rot_{\alpha_n}$ and  $\tilde h \circ h_n |\check \A \in U(h_n|\A, \alpha_n)$. In other words, $\tilde h|\check \A\in  \cU$ with:
\[ \cU:= U(h_n|\A, \alpha_n)\circ h_n^{-1}|\check \A\]
By  \cref{approxa}, the map $\tilde h$ is approximable modulo structural deformation. Hence 
there exists a diffeomorphism $h:  \A_{R} \hookrightarrow \A_\infty$ onto its image such that:
\begin{enumerate}[($M_1$)] \setcounter{enumi}{-1}
\item the map $h$ commutes with $\Rot_{\alpha_n}$. 
\item $h(\A)= \A$,   $ \det Dh|\A =1$ and  $h|\check \A\in \cU$, 
\item[($M_2'$)] the support of $ h$ is included in  $  \check \A_R= \check \A_\infty \cap \A_{R}$, 
\item[($M_3$)] the pullbacks  $h^* J_o $  and $h^*\Omega_o$ are  in  $\cW_J$ and $\cW_\Omega$,
\item[($M_4$)]  the maps $h$ and $\sigma$ commute.

\end{enumerate} 
Note that $h_{n+1}=h\circ h_n$ is well defined on $\A_{\rho}$ and supported by $\check \A_{\rho}$ by $(M_2')$ and $(V_n)$, as claimed by $(V_{n+1})$. 
 Also Property  $(III_{n+1})$ is satisfied by $(M_3)$. Property $(IV_{n+1})$ is an immediate consequence of $(M_4)$ and $(IV_n)$.  Note that for every $\alpha_{n+1}\in \Q/\Z$ close enough to $\alpha_n$, Property $(I_{n+1})$ is satisfied by $(M_1)$ and $(M_0)$.

Let us show $(II_{n+1})$.    
The map $h$ is defined on $\A_{R}$ whose interior contains both $h_n(\A_{\rho})$ and $\Rot_{\alpha_n} \circ h_n(\A_{\rho})$ by \cref{inclusion W 1,inclusion W 2}. 
By $(M_5)$, the map $h$ commutes with $ \Rot_{\alpha_n}$. So we have:
\[ \Rot_{\alpha_n} \circ h\circ h_n(\A_{\rho}) =    h\circ \Rot_{\alpha_n} \circ h_n(\A_{\rho})\; .\]
Consequently, for  $\alpha_{n+1}$ close to $\alpha_n$,  the compact set $ \Rot_{\alpha_{n+1}}\circ h\circ h_n(\A_{\rho'})$ is close to $h\circ \Rot_{\alpha_n}\circ h_n(\A_{\rho'})$. By $(II_n)$,  the compact set $h\circ \Rot_{\alpha_n}\circ h_n(\A_{\rho'})$  is included in the interior of $h\circ h_n(\A_{\rho})$. Consequently  for  $\alpha_{n+1}$ close to $\alpha_n$,  the compact set $ \Rot_{\alpha_{n+1}}\circ h\circ h_n(\A_{\rho'})$ is included in the interior of $h\circ h_n(\A_{\rho})$. In other words, the first part of $(II_{n+1})$ is satisfied.  Then  the map $h_{n+1} ^{-1} \circ \Rot_\alpha \circ h_{n+1}  $ depends continuously on  $\alpha$ close to $\alpha_n$, and equals $h_{n } ^{-1} \circ \Rot_{\alpha_n} \circ h_{n }  =f_n$ at $\alpha=\alpha_n$. From this we obtain the second part  of $(II_{n+1})$. 
\end{proof}
\subsection{AbC realizations on deformed real analytic symplectic spheres and disks}\label{AbC realization for deformed real analytic symplectic spheres and disks}\label{sec6.3}
  In this subsection, we prove   \cref{Main1} regarding the sphere and \cref{Main2}  in the disk case. They  state that given $\M\in \{\S, \D\}$,   every $C^\infty$-AbC scheme 
$(U,\nu)$ on $(\M,\Omega)$  is realizable by an analytic symplectomorphism for a  $C^\infty$-compatible real analytic symplectic structure.
These theorems are  consequences of the following proposition shown below.

\begin{proposition}\label{AbC_cyl_2}Let 
$\M\in \{\S, \D\}$ and  let $ \rho>\rho'>1$.  Then  for every $C^\infty$-AbC scheme $(U,\nu)$,  there is 
 a complexification $(\M_\rho, J)$ on which  $\Omega$ extends $J$-holomorphically   and  there is   a $J$-holomorphic  diffeomorphism  $f :   \M_{\rho'} \hookrightarrow \M_{\rho} $ which  extends a real symplectomorphism of $(\M,\Omega)$ realizing the   AbC scheme $(U,\nu)$.
 
 Moreover on the complement of $ \check \M_\infty$, the map $f$ equals  a rotation and   $J=J_o$. 
\end{proposition}

\begin{proof}[Proof of  \cref{Main1,Main2}   in the cases $\M=\S$ and  $\M=\D$]  
It is a direct consequence of  \cref{AbC_cyl_2} and respectively \cref{Def struture analytico symplectic} and  \cref{Def struture analytico symplectic2}. \end{proof}

\begin{proof}[Proof of \cref{AbC_cyl_1}] 
The space of smooth embeddings $\Emb^\infty (\M_{\rho}, \M_\infty)$ is    metrizable for a complete distance $d$.
The space  of complex differentiable  2-forms  on  $\M_{\rho}$  is  metrizable for a complete distance denoted also by $d$.
Likewise the space of smooth linear automophisms of $T^\R\M_{\rho}$ is    metrizable for a complete distance denoted again by~$d$.

We recall that $\check \M\approx \check \A$, hence a naive idea is to realize  the AbC scheme with maps $(h_n)_n$ which extend holomorphically   to a certain $\A_R$ -- as implemented in the previous section --  and then pull  back these extensions by $\pi_\M$. 
Yet there is an obstructions to this strategy. While the map $h_n \in \Emb^\infty(\A_\rho, \A_\infty)$ coincides with the identity on the complement of $\check \A_\rho$, the image $h_n(\check \A_\rho)$ can intersect the complement of $\check \A_\infty$, see \cref{sketch_Sphere} \cpageref{sketch_Sphere}. As the range of $\pi_\M$ is $\check \A_\infty$, this forbids us to lift $h_n$ by $\pi_\M$ as an embedding $\check \M_\rho\hookrightarrow \M_\infty$. Yet we remark that  $f_n=h_n^{-1} \circ \Rot_{\alpha_n} \circ h_n$ sends $\check \A_{\rho'}$ into $\check \A_\rho$, which enables us to pull back it by $\pi_\M$.   The pullback of $f_n$ by $\pi_\M$ is then defined on $\pi_\M^{-1} (\check \A_{\rho'})$ but the latter does not contain any $\M_{r}$, among $r>1$. However, for suitable $r$ and $\rho'$ as depicted \cref{sketch_Sphere}, we can extend smoothly  $f_n$ by $\Rot_{\alpha_n}$ on $\M_r\setminus \pi_\M^{-1} (\check \A_{\rho'})$. This is equivalent to say that $f_n$ equals $\Rot_{\alpha_n}$ outside of $\check \S_r$ and that $f_n|\check \S_r$ is semi-conjugate to $\Rot_{\alpha_n}$ by the map $H_n$ which equals $h_n\circ \pi_\M $ on $\check \S_r\cap \pi_\M^{-1} (\check \A_{\rho'})$ and to $\pi_\M$ on 
$\check \S_r\setminus \pi_\M^{-1} (\check \A_{\rho'})$. Note that $H_n$ is in general a non-injective local diffeomorphism, see \cref{sketch_Sphere}.

This  leads us to consider 
  the space  $\mathrm{Dilo}^\infty_c (\check\M_{\rho} , \A_\infty)$  of local diffeomorphisms $H: \check \M_{\rho}\to \A_\infty$ which coincide  with $\pi_\M$ on the complement of a compact subset of $\check \M_{\rho} $.  The map $H$   induces  the following  complex structure and 2-form on $\M_{\rho}$:
\[  H^\# J_o:= \left\{ \begin{array}{cl} 
  H^*  J_o  &\text{ on } \check \M_{\rho}\\
J_o &\text{ on } \M_{\rho}\setminus \check \M_{\rho}\end{array} \right.\qand   H^\# \Omega_o:= \left\{ \begin{array}{cl} 
 H^*\Omega_o &\text{ on } \check \M_{\rho}\\
\Omega_o &\text{ on } \M_{\rho}\setminus \check \M_{\rho}\end{array} \right.\; .\]
If furthermore 
$H(\check \M)=\check \A$ and 
$H|\check \M$ is symplectic, then we define  $\widehat H \in \Symp^\infty (\M)$   by:
\begin{equation} \label{hat}   \widehat H := \left\{ \begin{array}{cl} 
 \pi^{-1}_\M\circ H  &\text{ on } \check \M \\
id &\text{ on } \M \setminus \check \M \end{array} \right. \; .\end{equation} 
Let $\mathrm{Diro}^\infty_c (  \check \M_{\rho'},\M_{\rho})$ be the space of embeddings $f:  \M_{\rho'}\hookrightarrow \M_{\rho}$ which coincide  with a rotation  on the complement of a compact subset of $\check \M_{\rho'} $.  
Similarly to  \cref{example extension symp}, we have:
\begin{fact}\label{diese preservance} If $H\in \mathrm{Dilo}^\infty_c (\check\M_{\rho} , \A_\infty)$,
  then     $H^\#\Omega_o$ is $ H^\#J_o$-holomorphic. \end{fact} 

From now on, the proof is similar to that of \cref{AbC_cyl_1}; the difference is that $(h_n)_n$ is not formed by  diffeomorphisms from  $\A_{\rho} $ into $\A_\infty$  but \emph{local} diffeomorphisms $(H_n)_n$ from    $\check \M_{\rho}$ into $\A_\infty$.

We are going to construct by induction on $n\ge 0$ the existence  of $\alpha_n\in \Q/\Z$, $H_n\in \mathrm{Dilo}^\infty_c (\check\M_{\rho} , \A_\infty)$  and a diffeomorphism  $f_n\in  \mathrm{Diro}^\infty_c (  \check \M_{\rho'},\M_{\rho})$  such that:
\[ \alpha_0=0\; ,\quad  H_0= \pi |\A_{\rho}\quad\text{and so}\quad  f_0= id|\A_{\rho'}\; .\]
And for $n\ge 1$,  $\alpha_n$, $H_n$ and $f_n$ satisfy:
\begin{enumerate}[$(I_n)$]
\item   $H_n(\check \M)=\check \A$, the restriction  $H_n|\check \M$ is a symplectomorphism  and the map $\widehat H _n $ defined by \cref{hat}   respects the AbC scheme $(U,\nu)$  with $\alpha_n$:
\[\widehat H _n^{-1}\circ  \Rot_{\alpha_{n-1}}\circ \widehat h _n=\widehat H_{n-1}^{-1}\circ \Rot_{\alpha_{n-1}}\circ \widehat H_{n-1} \; ,\]\[ \quad 
\widehat H_{n}|\check \M   \in U(\widehat H_{n-1}, \alpha_{n-1})\qand |\alpha_{n }-\alpha_{n-1}|  < \nu(\widehat H_{n}, \alpha_{n-1})\; .
\] 
\medskip
\item The map $f_n$ sends $\check \M_{\rho'}$ into  the interior of $\check \M_{\rho}$, satisfies 
 $d(f_n,f_{n-1})<2^{-n}$
and 
 $ H_n\circ f_n= \Rot_{\alpha_n}\circ H_n$ on $ \check \M_{\rho'}$.
\medskip
\item  The complex structure $J_n:=H_n^\#  J_o  $ and the form $\Omega_n:= H_{n}^\#\Omega$ satisfy:
\[  d(J_n,J_{n-1})<2^{-n}\qand d(\Omega_n,\Omega_{n-1})<2^{-n}\; .\]
 \item  The maps $H_n$ and $\sigma$ commute:  $H_n\circ \sigma= \sigma\circ H_n$. 
\end{enumerate} 
 We will prove the following below:
 \begin{lemma}\label{induction2} If  the induction hypothesis $(I_n\cdots IV_n)$ is satisfied at step $n\ge0$, then it is satisfied at step $n+1$. 
\end{lemma} 
Now let us show that the induction hypothesis implies the proposition. 
By  $(II_n)$,  the sequence $(f_n)_n$ converges to a certain $f\in \Emb^\infty(\M_{\rho'}, \M_{\rho})$.  By $(I_n)$, the map $f|\M$ is a symplectomorphism of $\M$ which  realizes the AbC scheme.  

 By $(III_n)$,  the sequence $(J_n)_n$ converges to  an almost complex structure $J$. By \cref{cpx sr closed}, $J$ is a complex structure.  Also by $(I_n)$ and \cref{diese preservance}, the complex form $\Omega_n$ is a $J_n$-holomorphic extension of $\Omega$ to $\M_\rho$. 
By $(III_n)$ and  \cref{CS holo}, the sequence   $(\Omega_n)_n$ converges to  a $J$-holomorphic extension of $\Omega$. Also, it holds:
\[f_n^*J_n= (H_n^{-1} \circ \Rot_{\alpha_n}\circ H_n )^* H_n^*J_o= H_n^* \Rot_{\alpha_n}^*J_o=J_n\quad \text{ on } \check \M_\rho\; .\]
 Hence at the limit $f^*J=J$.  Consequently, the map $f$ is a  $J$-holomorphic extension  to $\M_{\rho'}$  of  the symplectomorphism  $f|\M$. 
  Finally as $H_n\in \mathrm{Dilo}^\infty_c (\check\M_{\rho} , \A_\infty)$  and  $f_n\in  \mathrm{Diro}^\infty_c (  \check \M_{\rho'},\M_{\rho})$, 
     it comes   that  on the complement of $ \check \M_\infty$, the map $f$ equals  a rotation and   $J=J_o$. 
\end{proof}
\begin{proof}[Proof of \cref{induction2}]
Let $\eta_o>0$ be so small that $\{ z\in \check \M_\rho: H_n(z)\neq \pi_\M(z)\}$ is disjoint from the $\eta_o$-neighborhood of $\M_\rho \setminus \check \M_\rho$. In particular:
\[ \{ z\in \check \M_\rho: H_n(z)\neq \pi_\M(z)\}\; \subset\;   \M_{{\rho}} \cap \pi_\M^{-1} \{ (\theta, y)\in \A_\infty : |\Re (y)|\le 1-\eta_o\}\; .\]

As the scheme $(U,\nu)$ is $C^\infty$-realizable, there exists $h_\M\in \Symp^\infty(\M)$ such that $  h_\M|\check \M$ 
commutes with $\Rot_{\alpha_n}$ and such that $  h_\M\circ \widehat H_n|\check \M$
belongs to $U(\widehat H_n, \alpha_n)$.
By density, we can assume $h_\M|\check \M$ compactly supported. Then there exists a compactly supported $h_\A\in \Symp^\infty(\check \A)$ such that $\pi_\M \circ   h_\M=h_\A\circ \pi_\M$ on $\check \M $. 
Then observe that $h_\A$ commutes with $\Rot_{\alpha_n}$ and that $  h_\M\circ \widehat H_n= \widehat{h_\A \circ H_n}$ on $\check \M$.  Then $h_\A|\check \A$ belongs to the following open set:
\[\cU:= \{ h \in \Symp^\infty (\check \A): \pi_\M^{-1} \circ h\circ H_n|\check \M \in U(\widehat H_n, \alpha_n)\}\; .\]

 By  \cref{approxa}, the map $h_\A$ is approximable modulo structural deformation. Hence 
there exists $\eta\in (0,\eta_o)$ such that  the following property holds true. 

For every $R>1 $, for any $C^\infty$-neighborhoods 
  $\cW_J$ of $J_o| \A_{R} $ and   $\cW_\Omega$ of $\Omega| \A_{R} $,  then  there exists a smooth diffeomorphism   $h:  \A_{R} \hookrightarrow \A_\infty$  onto its image  such that:
\begin{enumerate}[($M_1$)] \setcounter{enumi}{-1}
\item  the map $h$ commutes with $\Rot_{\alpha_n}$. 
\item $h(\A)= \A$,    $ \det Dh|\A =1$ and  $ h  |\check \A\in \cU$,  
\item the map  $ h$ is supported by   $ \{ (\theta, y)\in \A_{R}: |\Re (y)|<1-\eta\}$,
\item the pullbacks  $h^* J_o $  and $h^*\Omega_o$ are  in  $\cW_J$ and $\cW_\Omega$,
\item the map  $ h$    commutes with $\sigma$. 
\end{enumerate}
Let   $S_\rho:= \M_{\rho} \cap \pi_\M^{-1}\{ (\theta, y)\in \check \A_\infty: |\Re (y)|\le 1-\eta\}$. Note that $S_\rho$ is a compact subset of $\check \M_\rho$.  Hence we can fix    $R>1$ such that the interior of $\A_R$ contains  $ H_n(S_\rho)$:
\begin{equation} \label{inclusion W 1S}  \A_{R}\Supset  H_n(S_\rho) \; .\end{equation}
There is a $C^\infty$-neighborhood $  \cW_J$ of $J_o|  \A_{R}$  such that for every $J \in  \cW_J$ equal to $J_o$ on  $\{ (\theta, y): |\Re(y)|\ge 1-\eta\}$, the  structure $H_n^\#  J $ is   at distance $\le 2^{-n-1}$ from  $J_n =H_n^\#  J_o$. 
Also there is a $C^\infty $-neighborhood $ \cW_\Omega$ of $\Omega_o|  \A_{R} $ such that for every $\Omega'\in \cW_\Omega$ equal to $\Omega_o$ on $\{ (\theta, y): |\Re(y)|> 1-\eta\}$, the  form  $H_n^\# \Omega'$   is at distance $\le 2^{-n-1}$ from  $\Omega_n=H_n^\#  \Omega_o $. 

Let $h$ be satisfying $(M_0\cdots M_4)$ with these settings for   $\A_{R}$,  $ \cW_\Omega$ and $ \cW_J$.  Set:
\[H_{n+1} :(\theta,y)\in \check \M_\rho \mapsto \left\{ 
\begin{array}{cl} 
h\circ H_n(\theta,y) & \text{if } |\Re(y)| \le 1-\eta\; ,\\
\pi_\M(\theta,y) & \text{otherwise}\; .\end{array} \right. \]
By definitions of  $ \cW_\Omega$ and $ \cW_J$, the map   $H_{n+1}$ satisfies immediately  $(III_{n+1})$. By $(M_4)$ and $(IV_n)$, its satisfies  $(IV_{n+1})$. 
By $(M_0)$ and $(II_n)$, we have on $S_\rho$:
\[  \Rot_{\alpha_{n}}\circ H_{n+1} =\Rot_{\alpha_{n}}\circ h\circ  H_{n} =h\circ \Rot_{\alpha_{n}}\circ H_{n} = h \circ H_n\circ f_n= H_{n+1} \circ f_n\; .\]
On $\check \M_\rho\setminus S_\rho$, we have:
\[  \Rot_{\alpha_{n}}\circ H_{n+1} =\Rot_{\alpha_{n}}\circ \pi_\M =\pi_\M \circ   \Rot_{\alpha_{n}} =  H_{n+1} \circ f_n\; .\]
 Consequently:
 \begin{equation} \label{preII}
  \Rot_{\alpha_{n}}\circ H_{n+1}  =  H_{n+1} \circ f_n\quad \text{on } \check \M_\rho\; .
 \end{equation} 
By  $(M_1)$ and $(I_n)$, the restriction  $H_{n+1}|\check \M$ is a symplectomorphism onto $\check \A$. Furthermore, it satisfies   $(I_{n+1})$ with any $\alpha_{n+1}$ sufficiently close to $\alpha_n$. 

It remains to fix $\alpha_{n+1}$ close to $\alpha_n$ so that $f_{n+1}$ exists and satisfies $(II_{n+1})$.  
Let:
\[S_{\rho'}:= \M_{\rho'}\cap \pi_\M^{-1}(  \{ (\theta,y)\in \A_\infty:    |\Re (y)|\le1-\eta\}\; .\]
First recall that $f_n $ coincides with $\Rot_{\alpha_n}$ on the complement of  $S_{\rho'}$ and sends the  compact set $S_{\rho'}$ into the interior of $\check \M_{\rho}$ by $(II_n)$.  Thus, the map $H_{n+1}$ is defined on a neighborhood of $f_n(S_{\rho'})$. 
Hence for any small angle $\epsilon\in \T$, there exists a unique map $g_\epsilon: f_n(S_{\rho'})\to \M_\infty$ which is $C^\infty$-close to the canonical inclusion such that:
\[ H_{n+1} \circ g_\epsilon= \Rot_\epsilon \circ H_{n+1}\quad \text{on } f_n(S_{\rho'})\; .\]
By uniqueness and since $H_{n+1}=\pi_\M$ on $\check \M_{\rho'}\setminus S_{\rho'}=f_n(\check \M_{\rho'} \setminus S_{\rho'})$, we can extend smoothly $g_\epsilon$ by the rotation $\Rot_\epsilon$ on $f_n(\check \M_{\rho'}\setminus S_{\rho'})$. Then we have:
\[ H_{n+1} \circ g_\epsilon= \Rot_\epsilon \circ H_{n+1}\quad \text{on } f_n(\check \M_{\rho'} )\; .\]
Using the latter equality and \cref{preII}, for   $\alpha_{n+1}$ close to $\alpha_n$, with $\epsilon = \alpha_{n+1}-\alpha_n$ and $f_{n+1}:= g_\epsilon\circ f_n$, we have on $\check \M_{\rho'} $: 
\[\Rot_{\alpha_{n+1}} \circ H_{n+1} 
=\Rot_\epsilon \circ \Rot_{\alpha_n} \circ H_{n+1} 
=\Rot_\epsilon \circ H_{n+1}\circ f_n 
=H_{n+1} \circ g_\epsilon\circ f_n
= H_{n+1} \circ f_{n+1}\]
as sought by the last part of $(II_{n+1})$. The last part of $(II_{n+1})$ is then easily obtained by taking $\alpha_{n+1}$ close enough to $\alpha_n$.  
\end{proof} 

\subsection{Proof of the approximation theorem modulo structural deformation}\label{sec6.4}\label{proof approx}Let us show \cref{approxa}  which stated:

\begin{theorem*}
For every $h_o \in \Symp^\infty(\A)$,  for every  $C^\infty$-neighborhood $\cU\subset \Symp^\infty(\check \A)$  of $h_o|\check \A$, then there exists $\eta>0$ such that  the following property holds true. 

For every $R>1 $, for any $C^\infty$-neighborhoods 
  $\cW_J$ of $J_o| \A_{R} $ and   $\cW_\Omega$ of $\Omega_o| \A_{R} $,  then  there exists a smooth diffeomorphism   $h:  \A_{R} \hookrightarrow \A_\infty$  onto its image  such that:
 
\begin{enumerate}[($M_1$)]  
\item $h(\A)= \A$,    $ \det Dh|\A =1$ and  $h|\check \A\in \cU$,  
\item the map  $ h$ is supported by   $ \{ (\theta, y)\in \A_{R}: |\Re (y)|<1-\eta\}$,
\item the pullbacks  $h^* J_o $  and $h^*\Omega_o$ are  in  $\cW_J$ and $\cW_\Omega$,
\item  the map  $ h$  commutes with $\sigma$.
\end{enumerate}

\end{theorem*}

\label{proof of approxa}
To show this theorem, for $\rho>1$, we set:
\[K_\rho:= \T_\rho \times Q_\rho\quad \text{where }\T_\rho := \T+i [-\rho ,\rho]\qand Q_\rho:= [-\rho,-1  ]\sqcup [1  ,\rho] + i [-\rho ,\rho] 
\; .\]
We denote by  $\Ham^\infty_c (\check \A)$ the space of $C^\infty$-symplectomorphisms of $\check \A$ which are compactly supported  and isotopic to $id$ via a compactly supported  isotopy. 
Here is a the main approximation theorem of \cite{Be22}:
\begin{theorem} \label{approx analytic}Let  $g_o\in \Ham_c^\infty( \check \A)$   and  let    $\cU$ be  a $C^\infty (\check \A, \T\times \R)   $-neighborhood of  $g_o$.  
Then for any $\rho> 1$,  there exists an embedding   $ g_\rho|\check \A \in \cU $ which is the restriction of a biholomorphic map  $g_\rho$ of $\A_\infty$ leaving invariant $\Omega_o$, commuting with $\sigma$ and  such that:
\[ \sup_{(\theta,y)\in K_\rho} \| g_\rho(\theta,y)-(\theta,y)\|\le  \rho^{-1}\; .\]
 \end{theorem} 
 As $g_\rho$ commutes with $\sigma$, it leaves invariant $\T\times \R$. However this theorem does not imply \cref{approxa}  because in general  $g_\rho$ does not leave $\A$ invariant.

Given $\epsilon>0$ and $\rho>1$, let us denote $K_{\delta, \epsilon}$  the $\epsilon$-neighborhood of $K_\rho$. 
A consequence of the latter theorem  is:
 \begin{corollary} \label{coro approx analytic} Let  $h_o\in \Ham_c^\infty(\check \A)$   and  let    $\cV$ be  a $C^\infty (\check \A, \T\times \R)   $-neighborhood of  $h_o$.  
Then there is $\epsilon>0$ such that for any integer  $r\ge  1$,  there exists an embedding  $h_r|\check \A\in  \cV $ which is the restriction of a biholomorphic map  $h_r$ of $\A_\infty$ leaving invariant $\Omega_o$, commuting with $\sigma$  and such that:
\[   \| (h_r -id)|K_{r,\epsilon} \|_{C^r} \le  r^{-1}\; .\]
 \end{corollary} 
 \begin{proof}[Proof that \cref{approx analytic} implies 
 \cref{coro approx analytic}]
 By definition of the $C^\infty (\check \A, \T\times \R)   $-topology, there exists $s\ge 1$ such that   $\cV$ contains the following:
 \[ 
\{ h\in  C^\infty (\check \A, \T\times \R): \| D^k h- D^k h_o\| \le 1/s, \;  \forall k\le s,\; (\theta, y)\in  \T \times [-1+1/s, 1-1/s]   \}\; .\]
Let $\epsilon\in (0,\tfrac1{2s})$ be such that the support of $h_o$ is included in 
 $\T \times (-1+2\epsilon, 1-2\epsilon)$.  Let:
 \[H: (\theta, y)  \mapsto (\theta, (1-2\epsilon)^{-1} y) .\]  
Hence $H$ expands the $y$-coordinate, satisfies that $H(\T \times (-1+2\epsilon, 1-2\epsilon))=\check  \A$. Thus the following  is in $\Ham^\infty_c(\check \A)$:
 \[ g_0:= (\theta,y)\in \check \A \mapsto H \circ h_o\circ H^{-1}   (\theta,y) \]
   As   $H(\T \times[-1+1/s, 1-1/s])\Subset \check  \A$, the following is a $C^\infty$-neighborhood of $g_0$:
\[  \cU:= \{ 
H \circ h\circ H^{-1}  : h |\check \A\in \cV\}\; .\]
Also observe that $H$ sends $K_{r,2\epsilon}$ onto $K_{\rho'}$ with $\rho':= (r+2\epsilon)/(1-2\epsilon)$.  Hence 
by applying \cref{approx analytic} with $\rho>\rho'$, there exists an embedding   $ g_\rho|\check \A \in \cU $ which is the restriction of a biholomorphic map  $g_\rho$ of $\A_\infty$ leaving invariant $\Omega_o$, commuting with $\sigma$ and  such that:
\[ \sup_{(\theta,y)\in K_\rho} \| g_\rho(\theta,y)-(\theta,y)\|\le  \rho^{-1}\; .\]
Now put:
\[ h_r:= H^{-1}\circ g_\rho \circ H \; .\]
We notice that $h_r$ is a biholomorphic map of $\A_\infty$ satisfying that $h_r|\check \A$ is in $\cV$.  Also observe that:
\[ \sup_{ K_{r,2\epsilon} } \| h_r -id\| = \sup_{ H (K_{r,2\epsilon} ) } \| H^{-1} \circ g_\rho  -H^{-1} \|\le   \sup_{  K_\rho}  \|   g_\rho  -id \| \le  \rho^{-1}\; . \]
 Using the Cauchy inequality, for $\rho$-large enough, we obtain that $(h_r-id)| K_{r,\epsilon} $ is $1/r$-$C^r$-small. 
 \end{proof}
\begin{proof}[Proof of \cref{approxa}]
First note that being approximable modulo deformation is a $C^\infty$-closed property. As  $\Ham_c^\infty(\check \A)$ is $C ^\infty$-dense in $ \Symp^\infty(\check \A)$, we can assume that 
$h_o|\check \A\in \Ham_c^\infty(\check \A)$.  Let $\cU$ be a $C^\infty$-neighborhood of $h_o|\check \A$. Let $\cV\Subset \cU$ be a smaller $C^\infty$-neighborhood of $h_o|\check \A$.   By \cref{coro approx analytic},  there exist  $\epsilon>0$ and a sequence $(h_r)_{r\ge 1}$ of   biholomorphic maps  $h_r$ of $\A_\infty$ leaving invariant $\Omega_o$, commuting with $\sigma$, such that $h_r|\Check \A\in \cV$ and:
\begin{equation} \label{convergence vers 0}   \| (h_r -id)|K_{r,\epsilon}\|_{C^r}  \to 0\; .\end{equation} 
Now fix $\beta \in C^\infty(\R,\R)$ such that $\beta$ is supported by $(-1, 1) $ and $\beta=1$ on  a neighborhood  of $(-1+\epsilon, 1-\epsilon) $.   Let $R'>R$.  Let $r_0>1$ such that for every $r>r_0$, it holds 
$\A_{R'}\Subset \T_r \times ([-r, r]+i[-r,r])$ and so:
\[  \{ (\theta,y) \in \A_{R'}: |\Re(y)|\ge 1-\epsilon\}= \A_{R'}\cap  K_{r,\epsilon}\]

Using the Abelian group structure on $\A_\infty$, we define:
\[ \hat h_r: (\theta, y)\in \A_{R'}  \mapsto 
\beta\circ \Re(y)\cdot h_r(\theta,y) + (1-\beta\circ \Re(y) )\cdot (\theta,y)\in \A_\infty\; .\]
\begin{fact} 
When  $r$ large, the map $\hat h_r$ is a smooth diffeomorphism onto its image.
\end{fact} 
\begin{proof} 
By \cref{convergence vers 0},  $\hat h_r| \A_{R'}\cap  K_{r,\epsilon}$ is $C^\infty$-close to the canonical inclusion when $r$ is large.
On the other hand, $\hat h_r$ equals $ h_r$ on $\A_{R'}\setminus K_{r,\epsilon}$. Hence $\hat h_r$ is a local diffeomorphism when $r$ is large.  Let us prove the injectivity of $\hat h_r$, for $r$ large.  Let $\epsilon'<\epsilon$ be such that  $\beta=1$ on $[-1+\epsilon', 1-\epsilon']$. 

Let $U:= \A_{R'}\setminus K_{r,\epsilon'}$ and $V:= \A_{R'}\cap  K_{r,\epsilon}$. We notice that $\A_{R}= U\cup V$. Also $\hat h_r|U$ equals $ h_r$ which is injective, while $\hat h_r|V$ is $C^1$-close to the identity and so is injective. Hence to show the injectivity of $\hat h_r$,  it suffices to show that $\hat h_r(U\setminus V)$ is disjoint from $\hat h_r (V\setminus U)$. 

Note that $K_{r,\epsilon}$ is a neighborhood of the compact set  $V \setminus U$ because:
\[V \setminus U= (\A_{R'}\cap  K_{r,\epsilon} )\cap (\A_{R'}^c\cup  K_{r,\epsilon'})=   \A_{R'}\cap  K_{r,\epsilon'}\Subset K_{r,\epsilon}.\]
Then, as both $h_r|K_{r,\epsilon}$ and $\hat h_r|K_{r,\epsilon}$ are close to the canonical inclusion, when $r$ is large, it holds:
\[ \hat h_r(V\setminus U)\subset h_r(K_{r,\epsilon})\; .\]
As $h_r$ is a biholomorphism, the set  $\hat h_r(U\setminus V)=h_r(\A_{R'}\setminus K_{r,\epsilon}) 
$ is disjoint from $h_r(K_{r,\epsilon})$ and so from $\hat h_r(V\setminus U)$.  
\end{proof} 
Note also that $\hat h_r$   leaves $\A$ invariant. However, in general, we have 
 $ \det D\hat h_r|\A \neq 1$ (so $(M_1)$ is not satisfied).  Let $\eta>0$  be such that the support of $\beta $ is included in $(-1+\eta, 1-\eta)$.  Note that $\eta<\epsilon$. We have:
 \begin{fact} When $r$ is large enough, it holds:
 \begin{enumerate} 
 \item[$(M_2')$] the map  $\hat  h_r$ is supported by   $\A_{R'}\setminus K_{r, \eta }$ and  commutes with $\sigma$,
\item[$(M_3')$] $ \hat h_r^*J_o-J_o$ and $ \hat h_r^*\Omega_o-\Omega_o$ are $C^\infty$-small when $r$ large with support in 
 $K_{r, \epsilon }\setminus K_{r, \eta }$. \end{enumerate}
 \end{fact}
 \begin{proof} Property $(M_2')$ follows from the fact that $\beta$ is supported by $(-1+\eta, 1-\eta)$. 
Property  $(M'_3)$ follows from the fact that the restriction of $\hat h_r$ to $\A_{R'}\setminus K_{r, \epsilon}$ coincides with $h_r$, which  is holomorphic and leaves invariant $\Omega_o$ while the restriction of   $\hat h_r$ to $   \A_{R'}\cap K_{r,\epsilon}$ is close to the identity. 
\end{proof}

We show below:
\begin{lemma} \label{construction phi}
When $r$ is large, there exists $\phi_r:   \A_{R'}\to \A_\infty$ which is $C^\infty$-close to the canonical inclusion, which leaves  $\A$ invariant, which   commutes with $\sigma$, with support in $ \A_{R'}\setminus K_{r,\eta}\subset \check \A_\infty \cap \A_{R'}$
and such that:
\[ \phi_r^*\Omega_o|T\A= \hat h_r^*\Omega_o|T\A\; .\]
\end{lemma} 
Hence when $r$ is large, the range of $\phi_r$ contains $\A_{R}$.   This enables to define:
\[ h: (\theta,y)\in \A_{R}\mapsto  \hat h_r\circ   \phi_r^{-1}(\theta,y)\in \A_\infty\; .\]
Let us prove that $h$ verifies $(M_1)$.  We note that $h$  leaves invariant $\A$. Also on $T\A$: 
\[ h^*\Omega =  ( \hat h_r\circ   \phi_r^{-1})^* \Omega= (\phi_r^{-1})^* ( \hat h_r^*\Omega)=  (\phi_r^{-1})^*  (\phi_r^*\Omega)= \Omega\; .\]
Hence $\det Dh|T\A=1$. Also $h|\check \A $ is close to $h_o|\check \A$ when $r$ is large and so belong to  $\cU$. 

Let us show that $h$ verifies $(M_2)$. The supports of   both $\phi_r$ and  $\hat h_r$  are included in $\A_{R'}\setminus K_{r,\eta}$.  Hence the support of $h$ is included in:
\[   \A_{R} \setminus K_{r,\eta}= \{(\theta,y)\in \A_{R}:  |\Re(y)|\le 1-\eta\}\; .\]
Finally  $(M_3)$ follows from  $(M_3')$  and the fact that $\phi_r$ is close to the canonical inclusion.
We obtain $(M_4)$ by recalling  that   both $\hat h_r$ and $\phi_r$ commute with $\sigma$ and so $h$ does.  
\end{proof}
\begin{proof}[Proof of \cref{construction phi}]
Let $g\in C^\infty (\A, \R)$ be such that:
\[\hat h_r^*\Omega_o|T\A= g(\theta,y)\cdot d\theta\wedge dy\; .\]
We notice that $g$ is $C^\infty$ close to $1$ and  $g= 1$ on $\T\times \{ y \in [-1,1]: |y| > 1-\eta  \}$. We extend $g$ by $1$  on the remaining of  $\T \times \R$. Let:
\[\phi_1(\theta, y):= \left(\theta, \int_0^y g(\theta,t) dt-\tfrac12 \int_{-1} ^1 g(\theta,t) dt+1\right).\]
  We notice that $\phi_1$ is smooth  and satisfies $\det D\phi_1(\theta,y)= g(\theta,y)$. Hence on $T(\T\times \R)$:
\[ \phi_{1}^* \Omega_o|T(\T\times \A)  =  \det D\phi_1(\theta,y) \cdot d\theta\wedge dy= g(\theta,y)\cdot d\theta\wedge dy\; .\]
In particular $\phi_{1}^* \Omega_o|T\A = \hat h_r^*\Omega_o|T\A$.  We extend $\phi_1$ to $\A_\infty$ by:
\[\hat  \phi_1: (\theta, y)\in \A_\infty \mapsto  \phi_1( \Re(\theta), \Re(y)) +i \Im (\theta,y)\; .\]
The map $\hat \phi_1$ satisfies all the properties of the lemma but the fact that it may  not leave   $\A$ invariant and its support may not be included in $ \{ (\theta,y): |\Re y|\le 1-\eta\}$. 

Instead, there are two functions $\gamma_+$ and $\gamma_-$ in  $C^\infty(\T, \R)$ which are $C^\infty$-close to $1$ and $-1$ and such that for every $(\theta,y)\in \T\times \R$:
\[ 
\hat \phi_1(\theta,y) = (\theta, \gamma_+(\theta)+y)  \text{ if } y\ge 1-\eta\qand 
\hat \phi_1(\theta,y) = (\theta, \gamma_-(\theta)+y)  \text{ if } y\le -1+\eta\; .\]
We observe that:
\[ \int_\T \gamma_+ - \gamma_-  d\theta= \leb  \phi_1(\A) = \leb  \hat h_r(\A) =\leb \A=2\; .\]
Then there is $\tau_r\in \R$ small such that $\int_\T \gamma_+= 1+\tau_r$ and $\int_\T \gamma_- =-1+\tau_r$.  Let $\Gamma_\pm$ be a function of $\theta$ such that $\partial_\theta \Gamma_\pm = \gamma_\pm -\tau_r$.

Observe that $\Gamma_\pm$ are $C^\infty$-close to $\pm 1$ when $r$ is large. Hence using  a bump function, there is a small function $H_r: \T\times \R\to \R $ which is supported by $\T \times (\R\setminus [-1+3\eta, 1-3\eta])$  such that:
\[ H_r: (\theta, y) \in \T\times \R \mapsto \left\{ 
\begin{array}{cl}
 \Gamma_+(\theta)& \text{if } y\ge   1-2\eta \\
\Gamma_-(\theta)& \text{if } y\le   -1+2\eta  \end{array} \right.\quad . 
\] 
Let $\Fl_{H_r} ^1$ be the time one map of $H_r$ and observe 
that for every $(\theta,y)\in \T\times \R$:
\[ 
\hat \phi_1(\theta,y) = \Fl_{H_r} ^1 (\theta,y) +(0,\tau_r) \text{ if } y\ge 1-\eta\qand 
\hat \phi_1(\theta,y) = \Fl_{H_r} ^1 (\theta,y) +(0,\tau_r)  \text{ if } y\le -1+\eta\; .\]
Hence, with:
\[ \hat \phi_2:(\theta,y)\in \A_\infty\mapsto \Fl_{H_r} ^{-1}(\Re(\theta), \Re(y)-\tau_r) +i(\Im(\theta) , \Im(y))\; ,\]
the map $\phi_r:= \hat \phi_2\circ \hat \phi_1$ satisfies the requested properties.
\end{proof}
\begin{appendix} 
\section{Transitivity is AbC-realizable}\label{example of ABC}
In this appendix, we show that transitivity  is  an AbC realizable $C^0$-property.
Let us construct a $C^0$-AbC scheme which realizes a transitive symplectic map on the cylinder $\A$. The case $\M= \D$ or $\S$ are proved similarly.   We endow $\A$ with the distance:
\[ d((\theta, y) ,(\theta', y') )\le \max \{ |\theta-\theta'|, |y-y'|\}\]
so that its diameter is $1$.   We begin by establishing the following:
\begin{lemma}\label{preABC} For every $\epsilon>0$ and every integer $q\ge 1$, there exists  $h_{q,\epsilon}\in \Ham^\infty(\A)$ commuting with $\Rot_{1/q}$ and  such that 
$h_{q,\epsilon} (\T\times \{0\})$ is $\epsilon$-dense in $\check \A$. 
\end{lemma}
\begin{proof}
We first chose a Hamiltonian $H\in C^\infty(\A, \R)$  whose  time one map $\Fl^1_H$ sends $\T\times \{0\}$ to a curve which intersects both $\T\times \{1-\epsilon\}$ and $\T\times \{-1+\epsilon\}$. For instance take $H(\theta, y)= K\cdot  (1-y^2)\cdot \cos (2\pi \theta)$ for some large $K$.  
Then we take $M\ge 1/\epsilon$ and $\hat H(\theta,y):= H(q\cdot M\cdot \theta, y)$. Observe that the time one map $h_{q,\epsilon} $ of the flow of $\hat H$   satisfies the requested properties.
\end{proof}

Let us define the image of any $(h,\alpha)\in  \Symp^0(\A)\times \Q/\Z$ by our scheme. Let:
\[ f:= h^{-1}\circ  \Rot_{\alpha}\circ h\; .\]
 Let $j= j(h,\alpha)\ge 0$ be maximal such that   any point of $\A$ is at distance $\le 3^{-j}$ of a segment of $f$-orbit. Let    $N(h,\alpha)\ge 1$ be minimal such that 
  for a certain $x= x(h,\alpha)$,
 any point of   $ \A$ is at distance $\le 3^{-j}$ to $\{ f^k(x): 0\le k\le  N(h,\alpha)\}$.

 By uniform continuity, there is $\epsilon>0$ such that any pair of $\epsilon$-close points is sent by $h^{-1} $ to a pair of $3^{-j-2}$-close points.  Then $h^{-1} (C)$ is $3^{-j-2}$-dense for every $\epsilon$-dense curve $C\subset \check \A$.
Let $\alpha= p/q$. By  \cref{preABC}, the map:
\[\tilde  h = h_{q,\epsilon} ^{-1}\circ h  \]
satisfies
\begin{equation*} \tilde h ^{-1}\circ \Rot_{\alpha }\circ  \tilde h  =h ^{-1}\circ  \Rot_{\alpha }\circ h \; .\end{equation*} 
Also the curve $ \tilde h^{-1}(\T\times \{0\})$   is  $3^{- j-2}$-dense.  This property is actually $C^0$-open; there is a neighborhood $ U(h,\alpha)$ of $\tilde h|\check \A$, such that for every  $\hat h\in U(h,\alpha)$,  the curve $ \hat h^{-1}(\T\times \{0\})$   is  $3^{- j-2}$-dense. Let $\hat h\in U(h,\alpha)$ be  satisfying also:
\begin{equation}\label{equa premiere}   \hat h ^{-1}\circ \Rot_{\alpha }\circ  \hat h  =h ^{-1}\circ  \Rot_{\alpha }\circ h=f \; .\end{equation} 
There is $ \nu(\hat h, \alpha)$ such that    any   rational number $\hat  \alpha  \in (\alpha-\nu, \alpha +\nu)\setminus \{\alpha\}$ has its denominator so  large  that the map  $\hat f := \hat  h^{-1}\circ \Rot_{\hat \alpha}\circ \hat h$ displays a segment of orbit $(\hat f^k (\hat x ))_{0\le k\le \hat  N}$ whose  closed $3^{-j-1}$-neighborhood covers $\A$.  We assume $\hat N= N(\hat h, \hat \alpha)$ minima. Up to reducing $\nu(\hat h, \alpha)$,  we have that   $\hat N  \ge N(h,\alpha) $.  By   \cref{equa premiere}, we can assume $ \nu(\hat h, \alpha)$ even smaller such that each iterate $\hat f^k$ is $3^{-j-1}$-$C^0$-close to $f^k$, for every $-1\le k\le \hat N $.
\begin{proposition} 
The following is a $C^0$-AbC scheme:
\[ (h,\alpha) \in \Symp^0(\A)\times \Q/\Z \mapsto  (U(h,\alpha), \nu(h,\alpha) )\in   \cT\, ^0 \times (0,\infty)\; .\]
Moreover all the maps it constructs are transitive. 
\end{proposition} 
\begin{proof} 
Let  $(h_n)_n \in  \Symp^0(\A)^\N$ and $(\alpha_n)_n\in (\Q/\Z )^\N$  satisfying:
\[h_0=id\;,   \quad    \alpha_0=0\]
 and for $n\ge 0$:
\[h_n^{-1}\circ  \Rot_{\alpha_n}\circ h_n=h_{n+1}^{-1}\circ \Rot_{\alpha_n}\circ h_{n+1} \; ,\]\[ \quad 
h_{n+1}|\check \M   \in U(h_n, \alpha_n)\qand 0<|\alpha_{n+1}-\alpha_n|  < \nu(h_{n+1}, \alpha_n)\; .
\] 
Then by construction, the sequence $N_n:= N(h_n, \alpha_n)$  is  non-decreasing sequence and such that for every $n\ge 0$, the map $f_n$ displays a segment of orbit  $(f_n^k(x_n))_{0\le k\le  N_{n}}$  whose  closed $3^{-n}$-neighborhood covers $\A$. Also 
$ f_{n+1} ^k$ is $3^{-n-1}$-$C^0$-close to $f_n^k$, for every $-1\le k\le N_{n}$.

In particular $(f_n)_n$ is a Cauchy sequence and so converges to a certain $f\in \Symp^0(\A)$.   Moreover for every $n$, the iterates $(f^k)_{k\le N_{n}}$ is 
$  3^{-n} \le  \sum_{k> n}  3^{-k}$ close to $(f_{n}^k)_{k\le N_{n}}$. Thus it displays  a $  3^{-n+1} $-dense orbit. Hence $f$ is transtive.
\end{proof}
\section{Any AbC scheme can be forced to produce pseudo-rotations}
\label{section AbC pseudo} 
In this appendix, we show \cref{AbC pseudo} by proving that any $C^r$-scheme $(U,\nu)$ can be modified so that all its realizations are pseudo-rotations: symplectomorphisms with only 0, 1 or 2 periodic points when $\M=\A$, $\D$ or $\S$, and  which are moreover all elliptic. 

Observe that up to taking $\nu $ smaller, we can assume that for every realization
\[f= \lim_\infty h_n^{-1}\circ \Rot_{\alpha_n} \circ h_n, \]
with $\alpha_n = \tfrac {p_n} {q_n}$ and $p_n\wedge q_n=1$, then it holds $q_n\ge n$. 
Indeed it suffices to take $\nu(h, p/q)$ smaller than $\min\{ | \tfrac {p } {q}-  \frac\ell k | \neq 0:  {k\le q, \ell \in \Z} \}$ for every $p\wedge  q=1$.  
If we assume moreover  $\nu(h, p/q)< 2^{-q} \cdot q^{-q}$, then the limit $\alpha:= \lim \alpha_n$ is a Liouville number and thus irrational. 

Also observe that any map of the form $f_\bullet:= h^{-1}\circ \Rot_{p/q} \circ h$ does not have   periodic point  of period $< q$.

 When $\M=\A$, we put $\hat f_\bullet= f_\bullet$.  When $\M$ is $\S$ or $\D$, by \cref{Space}, there is  symplectomophism  $ \check \M\to \check \A$  which conjugates  $f_\bullet|\check \M$ to a symplectomorphism of $\check \A$  which extends to a $C^{r-1}$-symplectomophism $\hat f_\bullet $ of $\A= cl(\check \A)$. Indeed, the real trace of $\pi_\M$   extends to a symplectic blow-up $\A\to \M$  at the fixed points of $\Rot$. This implies that the following is positive:
\[ \epsilon(h, \tfrac p q):= 
\min\{ 
d(x, \hat f_\bullet^k(x)) :  x\in \A   \text{ and }    0< k<q
 \}\; .
\]
Then up to taking $\nu $ smaller, we can assume that  for every $\alpha$ which is $\nu (h, \tfrac pq)$-close to $\tfrac pq$, 
with $\hat f_{\bullet+1}:=  h^{-1}\circ \Rot_{\alpha} \circ h$,  it holds:
\[  \max_\A d(\hat f_\bullet^k(x), \hat f_{\bullet+1}^k (x))<  2^{-q}\cdot   \epsilon(h, \tfrac p q)\; , \quad \forall 0< k< q\; .\]
In particular, given a realization  $f=\lim f_n$ with $f_n:= h_n^{-1}\circ \Rot_{\alpha_n} \circ h_n$, applying the above inequality with $h=h_{n+1} $, 
$p/q= p_n/q_n$, 
$\alpha= \alpha_{n+1} $,   $f_\bullet := f_n$ and $f_{\bullet+1} := f_{n+1}$, it holds for every $0<k< q_n$:
\[  \max_\A d(\hat f_n^k(x), \hat f_{n+1} ^k(x)  )<  2^{-q_n} \epsilon(h_{n+1} , \tfrac {p_n}  {q_n} )\le  2^{-q_n} \max_\A  d(x, \hat f_n^k(x)    ) \; ,\]
\[  \max_\A  d(x, \hat f_{n+1}^k(x)    )\le \max_\A\left(  d(x, \hat f_n^k(x) )+d(\hat f_n^k(x), \hat f_{n+1} ^k(x)  ) \right) \le (1+2^{-q_n})  \max_\A  d(x, \hat f_n^k(x)   )\; .\]
Now we use inductively on $n> k$ the latter inequality and infer that $q_n\ge n$:
\[  \max_{\A}  d(x, 
  \hat f^k_{k+1} (x) 
  )\le     \prod_{n> k}  (1+2^{-q_n}) \max_\A d(\hat f^k(x), x) \le  \prod_{n> 0}  (1+2^{-n}) \max_\A d(x, \hat f^k(x) )  \; .\] 
As $ \prod_{n> 0}  (1+2^{-n}) <\infty$ and  $\max_{\A}  d(x, \hat 
  f^k_{k+1} (x) 
  )> 0$,  we obtain that $ \max_\A d(x, \hat f^k(x) )>0$ for every $k$.  From this we deduce that $f$ has only 0, 1 or 2 periodic points when $\M=\A$, $\D$ or $\S$. Also   these periodic points are fixed and with eigenvalues of modulus 1 and argument  $\alpha= \lim \alpha_n$ which is   irrational.  Hence the periodic points are all elliptic.

\end{appendix}
 
\bibliographystyle{alpha}
\bibliography{references.bib}

\end{document}